\documentclass[11pt]{article}

\usepackage[margin=1in]{geometry}
\usepackage{amsmath,amssymb,amsthm,mathtools}
\usepackage{enumitem}
\usepackage[hyphens]{url}
\usepackage{hyperref}
\usepackage{xurl}
\hypersetup{colorlinks=true,linkcolor=blue,citecolor=blue,urlcolor=blue,
  pdftitle={Wirsching's Positive-Predecessor-Density Program: Proofs of Conjectures 1 and 3},
  pdfauthor={Renato Augusto Tavares},
  pdfsubject={Number theory; 3n+1 dynamics; predecessor density},
  pdfkeywords={3x+1 problem; Collatz problem; predecessor density; positive asymptotic density; invariant density; two-scale difference equation; log-periodic asymptotics; uniform saddlepoint method; Mellin transform; Riemann zeta function}}
\usepackage{orcidlink}
\usepackage{fontawesome5}

\newtheorem{theorem}{Theorem}[section]
\newtheorem{proposition}[theorem]{Proposition}
\newtheorem{lemma}[theorem]{Lemma}
\newtheorem{corollary}[theorem]{Corollary}
\theoremstyle{definition}
\theoremstyle{remark}
\newtheorem{remark}[theorem]{Remark}
\newtheorem*{conjecture}{Wirsching's Conjecture 3}

\DeclareMathOperator{\artanh}{artanh}
\newcommand{\R}{\mathbb{R}}

\newcommand{\Z}{\mathbb{Z}}

\title{Wirsching's Positive-Predecessor-Density Program:\\
Proofs of Conjectures 1 and 3}

\author{Renato Augusto Tavares~\orcidlink{0009-0002-0196-3311}\\
\normalsize Universidade Federal de Goi\'as\\
\small \href{mailto:rat@discente.ufg.br}{rat@discente.ufg.br}\\
\small \href{https://orcid.org/0009-0002-0196-3311}{https://orcid.org/0009-0002-0196-3311}}
\date{\today}

\begin{document}

\maketitle

\begin{abstract}
Wirsching (2003) reduces uniform positive predecessor density for the
$3n+1$ map to a chain of five conditions, organized into three
conjectures. We prove two of them. Conjecture~1 concerns his
path-counting generators: in the convolution carrying them to his Elka
functions the partition weights grow subexponentially while the
binomial ratio decays geometrically, so a window of width
$O(\sqrt\ell)$ dominates and its radius fits inside the hypothesis.
Conjecture~3 concerns the asymptotics near $0$ of an invariant density
$\varphi$, a base-$3$ analogue of the Fabius density, against an
explicit $\varphi_0$ due to Berg and Kr\"uppel. The exact log-Laplace
transform of $\varphi$ splits into a smooth part, a $\log3$-periodic
correction $H$, and a doubly exponentially small remainder. Berg and
Kr\"uppel represented that correction as an infinite product in 1998;
their analysis did not determine whether it is constant. We give $H$ as a Fourier series with
coefficients in closed form in $\Gamma$ and $\zeta$; a classical
zero-free theorem for $\zeta$ shows it is not constant, and we enclose
its oscillation rigorously. Wirsching's comparison class fixes one
phase of $H$, and there Conjecture~3 holds with limit $e^{H(0)}$,
certified to lie in $(0.53412203666478,0.53412203666479)$. Off that class the phase sweeps a full
period, so the unrestricted asymptotic
$\varphi(t)\sim\kappa\varphi_0(t)$ fails. Condition $(\star4)$ follows, at every window radius, with $\mu=1/3$:
what Wirsching's argument needs is weaker than Conjecture~3 itself, and
the same saddlepoint chain settles it directly. With Conjecture~1 the
chain reduces to the single condition $(\star3)$. Conjecture~2 is his route to it and remains open.
\end{abstract}



\section{Introduction}
\label{sec:intro}

Let $T$ denote the $3n+1$ map on the positive integers, $T(n)=n/2$ for
$n$ even and $T(n)=(3n+1)/2$ for $n$ odd. Wirsching's
monograph~\cite{Wirsching1998Book} develops the predecessor-set
machinery this paper descends from, and~\cite{Wirsching2003} takes it
further, studying the density of predecessor sets through an averaging
construction on the $3$-adic integers. A family of operators
$S_\ell$, built from path counts in the predecessor graph, converges in
the strong operator topology to a limiting operator whose essential
part is a Markov chain on $[0,1]$. That chain has a unique invariant
density $\varphi$, a base-$3$ analogue of the Fabius density, and
Wirsching reduces uniform positive predecessor density to a chain of
five conditions $(\star1)$ through $(\star5)$ about $\varphi$ and the
path counts, organized into three conjectures. We are not aware of a
published resolution of any of the three.

Two of them are proved here.

Conjecture~1 asks that a pointwise lower bound on Wirsching's
generators, relative to their $3$-adic Haar average, transfer to his
Elka functions. Theorem~\ref{thm:wirsching-conj1} proves it. In the
convolution carrying the generators to the Elka functions the partition
weights grow subexponentially while the binomial ratio decays
geometrically, so a window of width $O(\sqrt\ell)$ dominates, and its
radius can be budgeted to sit inside the hypothesis. Together with
Wirsching's own two theorems this makes $(\star3)$ alone sufficient for
uniform positive predecessor density through proved implications only
(Corollary~\ref{cor:star3-density}).

Conjecture~3 concerns $\varphi$ near $0$, compared against an explicit
comparison function $\varphi_0$: the right-hand side of Berg and
Kr\"uppel's equation~(9.6) \cite{BergKruppel1998}, which they derive,
five years earlier, as the behaviour near $0$ of a solution of a
truncated version of the functional equation $\varphi$ satisfies.
Wirsching's equation~(7.11) defines $\varphi_0$ to be exactly that,
naming them in the sentence that introduces it, so the identification
is his and not a substitution made here. They gave the periodic
correction relating the two as an infinite product without ever
evaluating it. The exact log-Laplace transform of $\varphi$ decomposes
as a smooth part, a $\log3$-periodic correction $H$, and a doubly
exponentially small remainder (Theorem~\ref{thm:H}). We give $H$ as a
Fourier series with coefficients in closed form in $\Gamma$ and $\zeta$
(Proposition~\ref{prop:Hfourier}). The classical zero-free theorem for
$\zeta$ on $\operatorname{Re}s=1$ then shows $H$ is not constant, with
no computation, and its oscillation is enclosed rigorously
(Proposition~\ref{prop:Hcert}). Conjecture~3 is itself restricted: Wirsching asks for the limit
``uniformly for sequences $(z_\ell)\in\widetilde A_{\delta_5}$'', his
own comparison class, and that class fixes a single phase of $H$.
There the limit exists and equals $e^{H(0)}$, so \textbf{Conjecture~3
is true} (Theorem~\ref{thm:conj3}). Away from that class the phase
sweeps a full period, so the stronger unrestricted asymptotic
$\varphi(t)\sim\kappa\varphi_0(t)$, which Wirsching considered first
and declined, is false (Theorem~\ref{thm:unrestricted}). The two are
statements about the same ratio at two different scopes, and only the
restricted one is Conjecture~3. The function $e^H$ is Berg and
Kr\"uppel's own periodic factor for this eigenfunction
(Corollary~\ref{cor:bk}).

Condition $(\star4)$ follows, for every window radius, with $\mu=1/3$
(Theorem~\ref{thm:star4}). Wirsching's own route there is to compose
$(\star5)\Rightarrow(\star4)$ with Conjecture~3, and that works; what
his argument actually needs is weaker than Conjecture~3, and
Proposition~\ref{prop:wirschingreq} settles it directly from the same
saddlepoint chain, without the class bookkeeping. His chain is then
closed except at $(\star3)$, and Conjecture~2, the implication
$(\star4)\Rightarrow(\star3)$, is the only conjectural link left in it.

Conjecture~2 is not proved here and not attacked here. What can be said
about it from the outside is in \S\ref{sec:discussion}.

\S\ref{sec:related} places the problem against the density and
distributional literature. \S\ref{sec:chain} fixes Wirsching's objects
and states the chain from the primary source. \S\ref{sec:conj1} proves
Conjecture~1. \S\S\ref{sec:prelim} to \ref{sec:proofs} prove
Conjecture~3, through the periodic correction, a uniform saddlepoint
approximation, an envelope estimate, and the identification with Berg
and Kr\"uppel's asymptotic. \S\ref{sec:chain-closed} composes the
chain.

\section{Related work}
\label{sec:related}

The predecessor-density question is older than Wirsching's formulation
of it. Crandall raised the counting of $3n+1$ predecessors as a
quantitative problem \cite{Crandall1978}, and Sander improved his lower
bound on the number of descending integers below a given $x$
\cite{Sander1990}. Krasikov then turned the question into a system of
difference inequalities on the counting functions
\cite{Krasikov1989}. Applegate and Lagarias then produced density
exponents by two independent routes: a tree-search method following
Crandall \cite{ApplegateLagarias1995I}, and, separately, a
linear-programming solution of Krasikov's system
\cite{ApplegateLagarias1995II}. Krasikov and Lagarias later pushed the same
inequalities further, obtaining for each fixed $a\not\equiv0\bmod3$, and
for $x$ large depending on $a$, at least $x^{0.84}$ integers below $x$
with $a$ in their forward orbit \cite{KrasikovLagarias2003}. That is the
published record for the exponent. It is a sublinear count with
an $a$-dependent threshold, so it is weaker in both respects than the
uniform positive density, of order $x/a$, that Wirsching's chain
targets.

Wirsching's own contribution to this line begins with an improved
predecessor-set estimate \cite{Wirsching1993}, continues through the
combinatorial structure of those sets \cite{Wirsching1996}, the
monograph \cite{Wirsching1998Book} and the urn model
\cite{Wirsching1998Urns}, and produces in 2003 the five-condition chain
studied here \cite{Wirsching2003}.

Between the monograph and the chain he tried an analytic route.
Wirsching (2001) links an integral operator to the limiting behaviour of
transitions between counting functions attached to predecessor sets,
reduces the resulting integral equation to a linear retarded functional
differential equation, and reports that the normalized averages of the
Elka functions diverge pointwise to $+\infty$. He states there that what
the limiting behaviour of the transition means for the combinatorial
quantities themselves was not known \cite[\S\S1,5]{Wirsching2001}. The invariant density $\varphi$ that the
chain's last condition concerns belongs to a function class with its own
literature: Berg and Kr\"uppel supplied the asymptotic used here
\cite{BergKruppel1998}, and Rvachev's survey covers the
functional-differential family these densities sit in, which Wirsching's
2001 concluding remarks point at and leave open
\cite{Rvachev1990,Wirsching2001}. Operators attached to the Collatz map have also
been studied for their own dynamics, where the behaviour is opposite:
B\'ehani, building on earlier work, shows an associated operator on weighted
Bergman spaces is hypercyclic and, under conditions, chaotic
\cite{Behani2024}, while the averaging operator $W_3$ used below
contracts each affine hyperplane of fixed mass and fixes exactly one
function of unit mass supported in $[0,1]$.

A second line of work models $3n+1$ dynamics probabilistically, and it
is where $(\star3)$'s residue-level content lives. Borovkov and Pfeifer
estimated the Syracuse problem through a probabilistic model
\cite{BorovkovPfeifer2000}, Kontorovich and Sinai gave a structure
theorem for the maps in question \cite{KontorovichSinai2002}, and
Kontorovich and Lagarias set the stochastic models for $3x+1$ and
$5x+1$ side by side \cite{KontorovichLagarias2010}. Sinai studied the
inverse-path ensemble and its random-walk statistics
\cite{Sinai2003}, showed that the distribution his structure theorem
induces on $3^m$-progressions converges weakly to the uniform one
\cite[Theorem~1]{Sinai2003MMJ}, and proved the same for the residues
induced modulo $3^n$ by independent geometric variables, noting that
the strengthened version in which individual probabilities converge is
false \cite{Sinai2004}. Volkov's randomly labelled binary
tree for $5x+1$ carries direction-dependent edge laws
\cite{Volkov2006}, and Tao's Syracuse random variable \cite{Tao2022}
is a closely related canonical residue law, approaching uniformity at
fine scales while failing to be uniformly distributed on
$\mathbb Z/3^n\mathbb Z$. It is not the object $(\star3)$ constrains.
$(\star3)$ asks for a pointwise lower bound on fixed-cost slices of the
generators, at growing $3$-adic resolution and uniformly in
$a\in\mathbb Z_3^\times$; Tao's variable conditions on nothing of the
kind. The distance between the two is where the difficulty of
$(\star3)$ lives.

Current work runs in several directions. Computation has verified the
conjecture below $2^{71}$ \cite{Barina2025}. On the exponent,
$x^{0.84}$ \cite{KrasikovLagarias2003} remains the published record; a
recent preprint returns to the counting problem by a different route
and obtains $x^{0.3227}$, describing $0.84$ as the historical record
\cite{Liu2025}. A 2026 preprint of Mazur, unrefereed, raises the exponent to
$x^{0.90}$ by an exact feasible point for the level-$18$
Krasikov--Lagarias linear program together with an adaptive
elimination of advanced terms, with the whole theorem formalized in
Lean and its two native-computation steps replayed by independent exact
verifiers \cite{Mazur2026}. None of this line reaches what
Wirsching's chain targets: these are sublinear counts with an
$a$-dependent threshold, where the chain aims at a density of order
$x/a$ uniform in $a$. Other recent work addresses the finer structure of
trajectories, a different question again \cite{RozierTerracol2026}.

The third conjecture sits in a second literature, on the equation
$\varphi$ satisfies rather than on the $3n+1$ map. That equation is a
two-scale difference equation with dilation parameter $3$, in Berg and
Kr\"uppel's family~\cite{BergKruppel1998}, and $\varphi$ is a translate
of Rvachev's atomic function $h_3$~\cite{Rvachev1990}. Kato and
McLeod's study of $y'(x)=ay(\lambda x)+by(x)$~\cite{KatoMcLeod1971} is
the analytic ancestor. Wirsching noticed the connection to Rvachev's
functions himself and left it to future work~\cite{Wirsching2001}.

Solutions of such equations carry log-periodic fluctuations in their
asymptotics. The phenomenon goes back to de Bruijn on Mahler's
partition problem~\cite{DeBruijn1948} and to Erd\H{o}s and
Richmond~\cite{ErdosRichmond1976}, and the technique of extracting the
fluctuation's Fourier coefficients as residues of a Dirichlet series is
developed for radix-rational sequences by Dumas~\cite{Dumas2008} and
for $q$-regular sequences by Heuberger and
Krenn~\cite{HeubergerKrenn2018}. That work treats summatory functions
of sequences. Corollary~\ref{cor:bk} applies the same mechanism to an
eigenfunction of the difference equation itself, which is where the
$\Gamma$ and $\zeta$ coefficients of
Proposition~\ref{prop:Hfourier} come from, and with them the
non-constancy of $H$.

We are not aware of a published resolution of any of Wirsching's three
conjectures. The 2003
paper is cited sparsely, in bibliographies of the problem and in adjacent
work, and we found no published follow-up taking up any of the three.
For the
second, what this paper gives is an account of how far the machinery
natural to it reaches, and where the argument developed here stops.

\section{Wirsching's chain, from the primary source}
\label{sec:chain}

Every section, page, equation and result reference to
\cite{Wirsching2003} in what follows is to the published article, and
each was checked against it.\footnote{The numbering agrees throughout with the author's
preliminary version, circulated as \emph{On positive predecessor density
in $3n+1$ dynamics} (Katholische Universit\"at Eichst\"att), which we
also read. Sections~1 to~7 of the published article begin on pp.~773,
776, 778, 780, 780, 783 and~784; Definition~1 is on p.~772, Theorem~1 on
p.~775, Conjecture~1 on p.~778, Theorem~2 on p.~780, Theorem~3 on
p.~781, Theorem~6 and Corollaries~7 and~8 on p.~783, Conjecture~2 on
p.~784 and Conjecture~3 on p.~787. The item-by-item concordance is in the repository of the Data
Availability Statement.}

Wirsching's central combinatorial objects are the ``Elka'' functions
$e_\ell(k,a):=|E_{\ell,k}(a)|$, on domain
$\mathbb N_0\times\mathbb Z_3^\times$ \cite[\S1--2]{Wirsching2003},
where $E_{\ell,k}(a)$ is the set of
paths $b\xrightarrow{T}\cdots\xrightarrow{T}a$ in the Collatz graph
with $T=T_0$ exactly $k$ times and $T=T_1$ exactly $\ell$ times
($T_0,T_1$ the two branches of the
accelerated $3n+1$ map). Their normalized generators $\widetilde
g_\ell$ (defined below) live on the continuous state space
$\mathbb X=\mathbb I_3\times\mathbb Z_3^\times$: the
$\mathbb Z_3^\times$ factor carries the $3$-adic variable habitat of
the Syracuse measure, while the $\mathbb I_3\cong\{0,1,2\}^{\mathbb
N}$ factor carries the normalized step count $k/3^\ell$, converging to
a density $\varphi$.

The compactly supported \emph{generators} $g_\ell(k,a)$ are defined by
the recursion \cite[\S2, eq.~(2.1)]{Wirsching2003}
\begin{equation}\label{eq:g-recursion}
 g_{\ell+1}(k,a)=\sum_{0\le j<2\cdot3^\ell}
 g_\ell\bigl(k-j,\tfrac{2^{j+1}a-1}{3}\bigr),
 \qquad g_0(k,a)=\mathbf 1[k=0]\,\mathbf1[a\in\mathbb Z_3],
\end{equation}
linked to the Elka functions by $e_\ell=p_\ell*g_\ell$, where
$p_\ell(m)$ counts the ways to pay $m$ with coins $c_0=1$,
$c_j=2\cdot3^{j-1}$ ($j\ge1$) \cite[\S2]{Wirsching2003}. The class
$A_\delta$ consists of the integer sequences $(k_\ell)$ with
$|\ell-k_\ell|\le\delta\sqrt\ell$ for every $\ell$
\cite[eq.~(1.5)]{Wirsching2003}; $\widetilde A_\delta$ is the
analogous class of sequences $(x_\ell)\in\mathbb I_3^{\mathbb N}$, those
with $|\lfloor3^\ell x_\ell\rfloor-\ell|\le\delta\sqrt\ell$
\cite[eq.~(3.2)]{Wirsching2003}, used in conditions $(\star4)$ and
$(\star5)$ below.

Two densities remain. $\varphi$ satisfies the
untruncated equation $\lambda\varphi'(t)=a(\varphi(at)-\varphi(at-a+1))$
for $a=3$, $\lambda=2/3$ \cite[eq.~(7.9)]{Wirsching2003}; the base
parameter $a$ of that equation and of Berg and Kr\"uppel's is unrelated
to the $3$-adic unit $a$ everywhere else in this paper.

$\varphi_0$ is the right-hand side of Berg and Kr\"uppel's
\cite[Proposition~9.1, eq.~(9.6)]{BergKruppel1998}, which they derive as
the asymptotic behaviour, as $t\to0^+$, of a particular solution of the
truncated equation $\lambda g'(t)=ag(at)$. Their derivation is a
saddlepoint calculation whose hypotheses they do not verify: ``since the
check of these assumptions is only a question of routine, we drop it and
restrict ourselves to the necessary formal calculations''
\cite[\S9]{BergKruppel1998}. Nothing below depends on that: we take the
closed form itself as the \emph{definition} of $\varphi_0$ and prove
every comparison we need for it. At $a=3$ and $\lambda=2/3$ it reads
\begin{equation}\label{eq:phi0}
 \ln\varphi_0(t)=
 \varepsilon\ln(2\beta)-\tfrac12\ln(2\pi)
 +\gamma\ln t+\delta_{\mathrm{BK}}\ln(-\ln t)
 -\beta\ln^2\!\Bigl(\frac{t}{-\ln t}\Bigr),
\end{equation}
with
$\alpha=\tfrac12-\log_32$,
$\beta=1/(2\ln3)$,
$\delta_{\mathrm{BK}}=\tfrac12+\alpha-2\beta\ln(2\beta)$,
$\gamma=-2\beta-\delta_{\mathrm{BK}}-\tfrac12$ and
$\varepsilon=\tfrac12+\alpha-\beta\ln(2\beta)$. Proposition~9.1 states
the last three; $\alpha$ and $\beta$ come from
\cite[eq.~(9.4)]{BergKruppel1998}, and the specialization
$\lambda=(a-1)/a$ that fixes $\alpha$ for this equation sits between
Propositions 9.2 and 9.3 \cite[p.~180]{BergKruppel1998}, where the
printed line carries $a$ in place of $\ln a$; the value used here is the
one (9.4) forces. Their $\delta$ is written $\delta_{\mathrm{BK}}$
throughout, since $\delta$ is the window radius everywhere else here.

Berg and Kr\"uppel expect only that $\varphi/\varphi_0$ stays bounded and bounded away
from zero, weaker than convergence
\cite[\S9]{BergKruppel1998}, and Wirsching restates the truncated
equation and this asymptotic at
\cite[eq.~(7.10)--(7.11)]{Wirsching2003}; $\varphi_0$ is used in
Conjecture~3 below in place of $\varphi$ itself, which has no known
closed form.

\begin{proposition}[Base-$3$ Fabius density]\label{prop:fabius}
$\varphi$ is the density of $X=\sum_{j\ge1}2U_j\cdot3^{-j}$, $U_j$
i.i.d.\ uniform on $[0,1]$ (an identification not in the source), the
base-$3$ analogue of the Fabius density \cite{Fabius1966} (the
base\nobreakdash-$2$ case satisfies
$h'(x)=4h(2x)$ on $[0,1/2]$; here $\varphi'(x)=
(9/2)\varphi(3x)$ on $[0,2/3]$ by \cite[eq.~(7.7)]{Wirsching2003},
the same shape with base $3$ in place of base $2$: coefficient
$4=2^2/(2-1)$ becomes $9/2=3^2/(3-1)$, support cutoff
$1/2=(2-1)/2$ becomes $2/3=(3-1)/3$).
Extend every function on $[0,1]$ by zero to $\mathbb R$ before applying
$W_3f(x):=\tfrac32\int_{3x-2}^{3x}f(t)\,dt$. Among the functions $f$
with $\operatorname{supp}f\subset[0,1]$ and $\int_0^1f=1$, $\varphi$ is the
unique fixed point of $W_3$, and it is $C^\infty$, piecewise polynomial
away from the standard Cantor set \cite[Corollary~7]{Wirsching2003}.
For $f_0\in L^1([0,1])$ with
$\int_0^1f_0=1$, the
iterates $f_{n+1}:=W_3f_n$ converge geometrically,
$\|f_n-\varphi\|_1\le2^{-n+1}\|f_1-f_0\|_1$
\cite[Corollary~8]{Wirsching2003}.
\end{proposition}

\begin{proof}
Only the first claim needs an argument; the rest is quoted. Write
$X=\sum_{j\ge1}2U_j3^{-j}$ and $X'=\sum_{j\ge1}2U_{j+1}3^{-j}$, so
$X=(2U_1+X')/3$ with $X'\overset{d}{=}X$ independent of $U_1$. Here
$2U_1/3$ is uniform on $[0,2/3]$ and independent of $X'/3$, so the sum
of the two is absolutely continuous. Write $f$ for its density.
Conditioning on $U_1$ and differentiating,
\[
 f(x)=\frac{d}{dx}\,\mathbb{E}_{U}\bigl[F_{X}(3x-2U)\bigr]
     =3\int_0^1f(3x-2u)\,du
     =\frac32\int_{3x-2}^{3x}f(t)\,dt
     =W_3f(x),
\]
the substitution $t=3x-2u$ giving $du=-dt/2$. Each summand lies in
$[0,2\cdot3^{-j}]$, so $X\in[0,2\sum_{j\ge1}3^{-j}]=[0,1]$ and
$\operatorname{supp}f\subset[0,1]$; also $f\in L^1$ with
$\int_0^1f=1$. So $f$ is a fixed point of $W_3$ meeting the two
side conditions, and $\varphi$ is the only one
\cite[Corollary~7]{Wirsching2003}. Hence $f=\varphi$.
\end{proof}

Wirsching's paper reduces uniform positive predecessor density (that
the sets of predecessors of $a$ have lower asymptotic density bounded
below by $c/a$, for a single constant $c>0$ uniform over admissible
targets $a$) to a chain of five conditions
\cite[Definition~1, \S\S1--7]{Wirsching2003}:
\[
\begin{array}{ll}
\text{Theorem 1 (proved in \cite{Wirsching2003})} & (\star1)\Rightarrow
  \text{Positive Density}\\
\text{Conjecture~1} & (\star2)\Rightarrow(\star1)\\
\text{Theorem 2 (stated in \cite{Wirsching2003}, proof below)} & (\star3)\Rightarrow(\star2)\\
\text{Conjecture~2} & (\star4)\Rightarrow(\star3)\\
\text{Conjecture~3} & (\star5);\ (\star5)\Rightarrow(\star4)\text{ argued in
  \cite[\S7]{Wirsching2003}}
\end{array}
\]
That is the chain as the source leaves it. Its status after
\S\ref{sec:conj1} and \S\ref{sec:proofs} is the subject of
\S\ref{sec:chain-closed}: $(\star5)$ and $(\star4)$ are established
outright, every implication from $(\star3)$ downward is a theorem, and
Conjecture~2 is the one implication still conjectural. Here $(\star1)$ is a $\liminf$ lower bound on Wirsching's Elka functions
relative to their Haar average, uniform over sequences
$(k_\ell)$ with $|\ell-k_\ell|\le\delta\sqrt\ell$; $(\star2)$ is the same
statement for his generators $g_\ell(k,a)$; $(\star3)$ upgrades the
$\liminf$ over $\ell$ to a bound holding for every $\ell\ge\ell_0$, an
index Wirsching's Theorem~2 prints rather than leaves implied
\cite[\S3, p.~780]{Wirsching2003}; whether it is global in $a$ or one
per unit is not settled by the text and nothing below turns on it. And
$(\star4)$ asserts that there is a constant $\mu>0$ with
\[
 \liminf_\ell\frac{(W_3^\ell\chi_1)(x_\ell)}{(W_3^\ell\chi_0)(x_\ell)}
 \ge\mu,\qquad\text{uniformly for }(x_\ell)\in\widetilde A_\delta,
\]
with $\chi_0=\mathbf 1_{[0,2/3]}$, $\chi_1=\mathbf1_{[1/3,1]}$, and
$W_3$ the averaging operator of Proposition~\ref{prop:fabius}.

First, $(\star4)$ is a statement purely about the one-dimensional operator
$W_3$. It never mentions Wirsching's generators and never mentions the
$3$-adic unit $a$, so it supplies no residue-level information, which
is what $(\star3)$ is a statement about. The
source offers two bridges: the averaged-generator
identity $\lim_\ell\int_{\mathbb Z_3^\times}\widetilde g_\ell(x,a)\,da=
\varphi(x)$ for fixed $x\in\mathbb I_3$ \cite[\S6]{Wirsching2003}, and
Wirsching's Theorem~3 \cite[\S5]{Wirsching2003}, strong convergence
$S_\ell\to S_\infty$ for every $f\in C(\mathbb X)$, of the
transition operators $S_\ell:C(\mathbb X)\to C(\mathbb X)$ defined by
$\widetilde g_\ell=S_\ell(\widetilde g_{\ell-1})$
\cite[\S4]{Wirsching2003}, uniform on
bounded equicontinuous families. Neither bridge is uniform in growing
resolution, so neither closes the gap between $(\star4)$ and $(\star3)$.

Second, $(\star3)$ and $(\star2)$ are the same inequality, and the
implication between them is immediate. Wirsching states Theorem~2
without proof \cite[\S3]{Wirsching2003}, so we give one. Writing
$\widetilde g_\ell(x,a)=\gamma_\ell\,g_\ell(\lfloor3^\ell
x\rfloor,a)$ for a constant $\gamma_\ell>0$ not depending on $a$
\cite[\S3]{Wirsching2003}, the constant cancels from both sides of the
ratio each condition bounds, so $(\star3)$ reads
$g_\ell(\lfloor3^\ell x_\ell\rfloor,a)\ge\mu\,\bar g_\ell(\lfloor3^\ell
x_\ell\rfloor)$ for every $\ell\ge\ell_0$ and every
$(x_\ell)\in\widetilde A_\delta$. Given $(k_\ell)\in A_\delta$, put $x_\ell:=k_\ell3^{-\ell}$ wherever
that is legitimate. It is not legitimate at every level: $A_\delta$
admits $k_\ell$ outside $[0,3^\ell-1]$ at finitely many early levels,
where $k_\ell3^{-\ell}\notin\mathbb I_3$. Pick $L(\delta)$ beyond which
the whole window lies in $[0,3^\ell-1]$, which exists since the window
grows like $\sqrt\ell$ and the range like $3^\ell$. For $\ell\ge L$ set
$x_\ell:=k_\ell3^{-\ell}$, so $\lfloor3^\ell x_\ell\rfloor=k_\ell$ and
$(x_\ell)\in\widetilde A_\delta$ by \cite[eq.~(3.2)]{Wirsching2003}; for
$\ell<L$ take any admissible value. Then $(\star3)$ gives the bound at
every $k_\ell$ in the window for $\ell\ge\max(L,\ell_0)$, and a bound
holding from some index on has $\liminf$ at least as large, which is
$(\star2)$ with $\delta_1=\delta$ and $\mu_1=\mu$. The finite prefix is
invisible to a $\liminf$. What $(\star3)$ supplies over $(\star2)$ is
the stronger quantifier, a bound at every $\ell\ge\ell_0$ rather than
in the limit, and the implication just proved discards it.

Third, Theorem~1's conclusion, uniform positive predecessor density, is
about $a\in\mathbb N$, non-cyclic, $a\not\equiv0\pmod3$; its hypothesis
$(\star1)$, like $(\star2)$ and $(\star3)$, is stated over
$\mathbb Z_3^\times$.

\section{Wirsching's Conjecture~1}
\label{sec:conj1}

Throughout this section, fix $a\in\mathbb Z_3^\times$ and write
$\bar g_\ell(k):=\int_{\mathbb Z_3^\times}g_\ell(k,b)\,db$ and
$\bar e_\ell(k):=\int_{\mathbb Z_3^\times}e_\ell(k,b)\,db$ for the
$3$-adic averages. Put $c_0=1$ and $c_j=2\cdot3^{j-1}$ for $j\ge1$.
Wirsching's convolution formulas
\begin{equation}\label{eq:elka-convolution}
 e_\ell(k,a)=\sum_{j=0}^k p_\ell(k-j)g_\ell(j,a),
 \qquad \bar e_\ell=p_\ell*\bar g_\ell
\end{equation}
hold at each fixed $a$ and on the averages respectively
\cite[eq.~(2.3),~(2.4)]{Wirsching2003}, with $p_\ell(m)$ the number of
partitions of $m$ into the parts $c_0,\ldots,c_\ell$ with repetition.
Both $g_\ell$ and $p_\ell$ are non-negative.

The quantifier order matters here, and the source does not settle it.
Wirsching states each of $(\star1)$, $(\star2)$ and $(\star3)$ as
``if there are real numbers $\delta,\mu>0$ such that\ldots'', with $a$
occurring free inside the displayed inequality
\cite[\S\S1--3]{Wirsching2003}. He never writes the quantifier on $a$,
so nothing printed distinguishes
$\exists\delta_1,\mu_1\,\forall a$ from
$\forall a\,\exists\delta_1(a),\mu_1(a)$. We adopt the first, and the
reason is his own Definition~1: uniform positive predecessor density
there means one constant $c>0$ with
$\liminf_n|\{x\in P_T(a):x\le n\}|/n\ge c/a$ \emph{for each} $a$, the
constant standing outside the quantifier on $a$
\cite[\S1]{Wirsching2003}. A $\mu_1$ free to depend on $a$ and tend to
$0$ could not deliver that single $c$ through his Theorem~1, so the
global reading is the only one under which the chain does what he says
it does. Stated in full, and understood as a reading rather than as
something the source prints:
\[
 \exists\,\delta_1,\mu_1>0\ \ \forall a\in\mathbb Z_3^\times\ \
 \forall\varepsilon>0\ \ \exists\,\ell_0(a,\varepsilon)\ \
 \forall\ell\ge\ell_0\ \ \forall(k_m)\in A_{\delta_1}:\quad
 \frac{g_\ell(k_\ell,a)}{\bar g_\ell(k_\ell)}\ \ge\ \mu_1-\varepsilon .
\]
The stated uniformity is over sequences, so the eventual index
$\ell_0$ may depend on $a$; the constant $\mu_1$ may not. Nothing in
\S\ref{sec:conj1} turns on this: Theorem~\ref{thm:wirsching-conj1} is
stated at a fixed $a$ and carries the same $\mu_1$ from the generators
to the Elka functions there, so it holds under either reading, and the
discussion above only says which reading makes the chain deliver
Wirsching's Theorem~1.

\begin{lemma}[Window form of the hypothesis]\label{lem:window-form}
Suppose $(\star2)$ holds with radius $\delta_1$ and constant $\mu_1$.
Then for every $\varepsilon>0$ there is $\ell_0$ such that for every
$\ell\ge\ell_0$ and every integer $j$ with
$|j-\ell|\le\delta_1\sqrt\ell$,
\[
 g_\ell(j,a)\ \ge\ (\mu_1-\varepsilon)\,\bar g_\ell(j).
\]
\end{lemma}

\begin{proof}
$A_{\delta_1}$ is a product: a sequence lies in it exactly when
$|m-k_m|\le\delta_1\sqrt m$ at each $m$ separately
\cite[eq.~(1.5)]{Wirsching2003}, with no constraint tying different
levels. Write $W_m$ for the window of admissible indices at level $m$, a
finite set. On $0\le j\le3^m-m-1$ the average $\bar g_m(j)$ counts urn
occupancies summing to $j$ and is positive, but $W_m$ need not sit
inside that range at small $m$: at $\delta_1=1$ and $m=1$ the window is
$\{0,1,2\}$ while the support stops at $1$, and the ratio at $j=2$ is
$0/0$. Since $W_m$ grows like $\sqrt m$ and the support like $3^m$,
there is $L(\delta_1)$ with $W_m\subset[0,3^m-m-1]$ for every $m\ge L$.
Put
\[
 k^*_m:=\operatorname*{arg\,min}_{j\in W_m}
 \frac{g_m(j,a)}{\bar g_m(j)}\quad(m\ge L),
 \qquad k^*_m:=m\quad(m<L),
\]
which is well defined, and $(k^*_m)\in A_{\delta_1}$. Applying the
hypothesis to this one sequence,
\[
 \liminf_m \frac{g_m(k^*_m,a)}{\bar g_m(k^*_m)}\ \ge\ \mu_1 ,
\]
so there is $\ell_0$ with
$g_\ell(k^*_\ell,a)\ge(\mu_1-\varepsilon)\bar g_\ell(k^*_\ell)$ for
$\ell\ge\ell_0$; enlarge $\ell_0$ to $\max(\ell_0,L)$ so that
$k^*_\ell$ is the window argmin there. By the choice of $k^*_\ell$,
every $j\in W_\ell$ has
\[
 \frac{g_\ell(j,a)}{\bar g_\ell(j)}\ \ge\
 \frac{g_\ell(k^*_\ell,a)}{\bar g_\ell(k^*_\ell)},
\]
which is the claim.
\end{proof}

The lemma turns control of one index per level into control of the
whole window at one level, which is what the convolution needs. It uses
the hypothesis in its weakest reading, one sequence at a time, asking
only that the $\liminf$ along each admissible sequence be at least the
same $\mu_1$; so Theorem~\ref{thm:wirsching-conj1} does not depend on
how ``uniformly for sequences'' is construed.

\begin{lemma}[Tail of the convolution]\label{lem:conv-tail}
Fix $\delta>0$. There are absolute $\rho\in(0,1)$ and $A>0$ such that
for every $\tau>0$ and all $\ell$ beyond a threshold depending only on
$\delta$,
\[
 \frac{1}{\bar e_\ell(k)}\sum_{m\ge\tau\sqrt\ell}p_\ell(m)\,
 \bar g_\ell(k-m)\ \le\ \frac{A}{1-\rho}\,\rho^{\,\tau\sqrt\ell}
 \qquad\text{uniformly for }|k-\ell|\le\delta\sqrt\ell .
\]
In particular the left side is $O\bigl(e^{-c\sqrt\ell}\bigr)$ with
$c=\tau\log(1/\rho)>0$.
\end{lemma}

\begin{proof}
Cancelling $\prod_{j=0}^\ell(1-z^{c_j})$ against its inverse in
\[
 Q_\ell(z):=\sum_{k\ge0}q_\ell(k)z^k
 =\prod_{j=0}^\ell\frac{1-z^{c_j}}{1-z},\qquad
 P_\ell(z):=\sum_{m\ge0}p_\ell(m)z^m
 =\prod_{j=0}^\ell(1-z^{c_j})^{-1}
\]
gives $P_\ell(z)Q_\ell(z)=(1-z)^{-\ell-1}$, so by \eqref{eq:elka-convolution}
\begin{equation}\label{eq:ebar-binomial}
 \bar e_\ell(k)=\frac{1}{2\cdot3^{\ell-1}}\binom{k+\ell}{\ell},
\end{equation}
where $q_\ell(k)=2\cdot3^{\ell-1}\bar g_\ell(k)$ counts the
distributions of $k$ indistinguishable balls into urns
$U_0,\ldots,U_\ell$ with occupancy range $0,\ldots,c_j-1$
\cite[\S2]{Wirsching2003}, \cite[p.~980]{Wirsching1998Urns}, which is
what the support bound $\sum_j(c_j-1)=3^\ell-\ell-1$ forces. The urn
source writes the maximum occupancy as $2\cdot3^{j-1}-1=c_j-1$, so
$c_j$ here is that maximum plus one and not the capacity in its sense.

Terms with $m>k$ contribute nothing: $g_\ell(j,\cdot)=0$ for $j<0$, so
$\bar g_\ell(k-m)=0$ there, and the sum may be read over
$0\le m\le k$. Since $p_\ell(0)=1$ and every term is non-negative,
$\bar g_\ell(j)\le\bar e_\ell(j)$ for every $j$, so by
\eqref{eq:ebar-binomial}
\[
 \frac{\bar g_\ell(k-m)}{\bar e_\ell(k)}
 \le\frac{\binom{k-m+\ell}{\ell}}{\binom{k+\ell}{\ell}}
 =\prod_{i=1}^{m}\frac{k-i+1}{k+\ell-i+1}
 \le\Bigl(\frac{k}{k+\ell}\Bigr)^{m}.
\]
In the window $k\le\ell+\delta\sqrt\ell$, so
$k/(k+\ell)\le(\ell+\delta\sqrt\ell)/(2\ell+\delta\sqrt\ell)\to\tfrac12$
and the ratio is at most $\rho_0:=3/5$ for all large $\ell$, uniformly
in $k$. Fix $\rho:=7/10$. The partition count obeys
for an absolute constant $C$, independent of $\ell$ and $m$,
\begin{equation}\label{eq:partition-bound}
 p_\ell(m)\le\exp\bigl(C\log^2(m+2)\bigr):
\end{equation}
$p_\ell(0)=1$, and for $m\ge1$ only the parts $c_j\le m$ contribute, at
most $\log_3(m/2)+2$ of them, each with multiplicity at most $m$, so
$p_\ell(m)\le(m+1)^{\log_3(m/2)+2}$. The bound is uniform in $\ell$.
Hence $A:=\sup_{m\ge0}\exp\bigl(C\log^2(m+2)\bigr)(\rho_0/\rho)^m$ is
finite, being the supremum of a subexponentially growing factor against
a geometric one with ratio $6/7<1$, and
\[
 \sum_{m\ge\tau\sqrt\ell}p_\ell(m)\frac{\bar g_\ell(k-m)}{\bar e_\ell(k)}
 \le\sum_{m\ge\tau\sqrt\ell}p_\ell(m)\rho_0^{\,m}
 \le A\sum_{m\ge\tau\sqrt\ell}\rho^{\,m}
 =\frac{A}{1-\rho}\,\rho^{\,\lceil\tau\sqrt\ell\rceil}. \qedhere
\]
\end{proof}

\begin{theorem}[Wirsching's Conjecture~1]\label{thm:wirsching-conj1}
In the notation of \S\ref{sec:chain}, fix $a\in\mathbb Z_3^\times$. If
$(\star2)$ holds at $a$ with window radius $\delta_1$ (Wirsching's own
naming, \cite[\S2]{Wirsching2003}) and constant $\mu_1$, then
$(\star1)$ holds at $a$ with the same constant $\mu_1$ and some smaller
radius $\delta<\delta_1$, which may be taken to be $\delta_1/3$ and
depends on $\delta_1$ alone.

In particular the implication holds under either reading of the
quantifier on $a$: if $\delta_1$ and $\mu_1$ in $(\star2)$ are global in
$a$, so are $\delta$ and $\mu_1$ in $(\star1)$; and if they are allowed
to depend on $a$, the implication still holds at each $a$ separately.
Under the global reading, which is the one Wirsching's Definition~1
requires, this is exactly Conjecture~1.
\end{theorem}

\begin{proof}
Choose $\tau>0$ and $\delta>0$ with $\delta+\tau<\delta_1$, which is
possible with $\delta=\tau=\delta_1/3$, and let
$\varepsilon\in(0,\mu_1)$. All statements below are for $\ell$ large
enough that $\tau\sqrt\ell\ge1$ and $\ell\ge(\delta+\tau)^2$, so that
the retained range is nonempty and $k-m\ge0$ throughout it. Take $\ell_0$ from Lemma~\ref{lem:window-form} for
that $\varepsilon$. Let $|k-\ell|\le\delta\sqrt\ell$ with
$\ell\ge\ell_0$ large enough for Lemma~\ref{lem:conv-tail}. For
$m<\tau\sqrt\ell$ the index $j=k-m$ satisfies
$|j-\ell|\le(\delta+\tau)\sqrt\ell<\delta_1\sqrt\ell$, so
Lemma~\ref{lem:window-form} applies at every such $j$ at once. Dropping
the non-negative terms with $m\ge\tau\sqrt\ell$ from
\eqref{eq:elka-convolution} and then restoring them inside the average,
\[
\begin{aligned}
 e_\ell(k,a)
 &\ \ge\ \sum_{m<\tau\sqrt\ell}p_\ell(m)\,g_\ell(k-m,a)
 \ \ge\ (\mu_1-\varepsilon)\sum_{m<\tau\sqrt\ell}p_\ell(m)\,
        \bar g_\ell(k-m)\\
 &\ =\ (\mu_1-\varepsilon)\Bigl(\bar e_\ell(k)
       -\sum_{m\ge\tau\sqrt\ell}p_\ell(m)\,\bar g_\ell(k-m)\Bigr)
 \ \ge\ (\mu_1-\varepsilon)\bigl(1-\eta_\ell\bigr)\,\bar e_\ell(k),
\end{aligned}
\]
with $\eta_\ell=O(e^{-c\sqrt\ell})$ by Lemma~\ref{lem:conv-tail},
uniformly over the window and independent of $a$, every quantity
entering it being a Haar average. So for every sequence
$(k_\ell)\in A_\delta$,
$\liminf_\ell e_\ell(k_\ell,a)/\bar e_\ell(k_\ell)\ge\mu_1-\varepsilon$,
uniformly over such sequences because $\ell_0$ and $\eta_\ell$ are.
Letting $\varepsilon\downarrow0$ gives $(\star1)$ at radius $\delta$
with constant $\mu_1$, and $(\star1)$ asks only for some positive
constant.
\end{proof}

\begin{corollary}\label{cor:star3-density}
Condition $(\star3)$ implies uniform positive predecessor density on the
non-cyclic integers $a\not\equiv0\bmod3$.
\end{corollary}

\begin{proof}
Wirsching's Theorem 2 takes $(\star3)$ to $(\star2)$, stated in
\cite[\S3]{Wirsching2003} and proved in \S\ref{sec:chain} above;
Theorem~\ref{thm:wirsching-conj1} takes
$(\star2)$ to $(\star1)$, and his Theorem 1 takes $(\star1)$ to the
conclusion \cite[\S1]{Wirsching2003}. Both of his theorems quantify
their window radius and constant existentially, so the possibly smaller
radius produced by Theorem~\ref{thm:wirsching-conj1} is admissible in
Theorem 1 as it stands.
\end{proof}

Before this the chain from $(\star3)$ downward passed through
Conjecture~1. It no longer does, and $(\star4)\Rightarrow(\star3)$ is
the only conjectural link left below $(\star4)$.

\begin{remark}[Relation to the source]\label{rem:wirsching-conj1-source}
Wirsching records the partition-counting identity
$\prod_{j=0}^\ell(1-z^{c_j})^{-1}=\sum_mp_\ell(m)z^m$
\cite[\S2, eq.~(2.2)]{Wirsching2003} and the closed
form for $\bar e_\ell(k)$ used above \cite[\S1, eq.~(1.4)]{Wirsching2003},
and describes, without a generating function, the urn count
$q_\ell(k)=2\cdot3^{\ell-1}\bar g_\ell(k)$, crediting the combinatorial
identification to \cite{Wirsching1998Urns} \cite[\S2]{Wirsching2003};
the generator apparatus itself is developed in
\cite[Ch.~IV, p.~103 and Thm.~IV.1.14]{Wirsching1998Book}, the pinpoints
Wirsching himself gives for it \cite[\S2]{Wirsching2003};
the generating function
$\sum_kq_\ell(k)z^k=\prod_{j=0}^\ell(1-z^{c_j})/(1-z)$ used above does
not appear there. Wirsching states the implication itself as
Conjecture~1 \cite[\S\S1--2]{Wirsching2003}. The cancellation
$P_\ell(z)Q_\ell(z)=(1-z)^{-\ell-1}$ is what makes $\bar e_\ell(k)$
binomial.
\end{remark}

The second conjecture is untouched by this. Its antecedent
$(\star4)$ comes from Conjecture~3, an asymptotic assertion on
$\varphi(z_\ell)/\varphi_0(z_\ell)$ along the central-limit window
$|\ell-k_\ell|\le\delta_5\sqrt\ell$, which \S\S\ref{sec:prelim} to
\ref{sec:proofs} prove and \S\ref{sec:chain-closed} composes with the
present section.

\section{The invariant density and its Laplace transform}
\label{sec:prelim}

Proposition~\ref{prop:fabius} identifies $\varphi$ as the density of
$X=\sum_{j\ge1}2\cdot3^{-j}U_j$ with $U_j$ i.i.d.\ uniform on $[0,1]$,
the unique fixed point of $W_3$ of unit mass supported in $[0,1]$. At
dilation parameter $a=3$ and $\lambda=2/3$ it is an instance of Berg
and Kr\"uppel's family $\lambda\,\psi'(t)=a\bigl(\psi(at)-\psi(at-a+1)\bigr)$,
and a translate of Rvachev's atomic function $h_3$ \cite{Rvachev1990}.
Three further facts are used throughout below.

The reflection $x\mapsto1-x$ relating $\varphi$ to
\cite[Cor.~7]{Wirsching2003}'s fixed point is vacuous:
$\varphi(x)=\varphi(1-x)$ identically, since $X$ and
$1-X=\sum_{j\ge1}\tfrac{2}{3^j}(1-U_j)$ are equidistributed, because
$1-U_j\sim U_j$ independently for every $j$ and
$\sum_{j\ge1}2\cdot3^{-j}=1$. It
is continuous, by Lemma~\ref{lem:tail} below. And it is constant, equal
to $3/2$, on the middle third $[1/3,2/3]$: there $3x-2\le0\le1\le3x$,
so $\operatorname{supp}\varphi\subset[0,1]$ collapses the integral in
$W_3\varphi=\varphi$ to $\int_0^1\varphi=1$.

For $s > 0$ write
\[
K(s) = \log \mathbb{E}\bigl[e^{-sX}\bigr] = \sum_{j \ge 1} g\!\left(\frac{2s}{3^j}\right), \qquad
g(b) = \log\frac{1-e^{-b}}{b}.
\]
The branch of $g$ is the one analytic on $\{\operatorname{Re} w > 0\}$ and real on $(0,\infty)$;
since $(1-e^{-w})/w$ is analytic and non-vanishing there (its zeros lie on $2\pi i \Z \setminus
\{0\}$, and $w=0$ is a removable singularity with value $1$), and the half-plane is simply
connected, $g$ exists and is unique. Differentiating termwise, justified because $g(w) = -w/2 +
O(w^2)$ as $w \to 0$ so the tail of the series is normally convergent together with every
derivative,
\[
t(s) := -K'(s), \qquad V(s) := s^2 K''(s) > 0,
\]
the second an exact positivity, $K''(s)$ being the variance of $X$ under its exponential tilt by
$e^{-sX}$. Since $t'(s) = -K''(s) < 0$, $t$ is strictly decreasing on $(0,\infty)$; at the
endpoints, $t(0^+) = \mathbb E[X] = \sum_{j\ge1} 2\cdot3^{-j}\cdot\tfrac12 = \tfrac12$, and
$t(s)\to0$ as $s\to\infty$ because the exponential tilt concentrates on $\operatorname{ess\,inf}
X=0$ (each factor $e^{-sX}$ suppresses any mass away from $0$ as $s\to\infty$, and $X\ge0$ a.s.).
So $t(s)$ decreases strictly from $1/2$ to $0$, and studying $\varphi(t(s))$ as $s \to \infty$ is
the same as studying $\varphi(t)$ as $t \to 0^+$.

The comparison function $\varphi_0$ is \eqref{eq:phi0}, the right-hand
side of Berg and Kr\"uppel's~(9.6), which they derive as the behaviour
near $0$ of a solution of the truncated equation
$\lambda\,\psi'(t)=a\,\psi(at)$, by Laplace transform and saddlepoint
analysis, with the constants $\alpha,\beta,\gamma,\delta_{\mathrm{BK}},
\varepsilon$ fixed there. Their dilation parameter $a$ is reused in
\S\ref{sec:periodic} for this paper's own, unrelated linear coefficient
of $Q$; the two never appear in the same formula except in
\S\ref{sec:bk}, where the collision is flagged at the point it matters.
Write $\varphi_{0,\mathrm{bare}}$ for \eqref{eq:phi0} without the
leading constant $(2\beta)^\varepsilon/\sqrt{2\pi}$.

\section{The periodic correction}
\label{sec:periodic}

Set $c = \log 3$ and write $w = \log s$, $L(w) = K(e^w)$. Substituting into the defining series and
splitting off the $j=1$ term against the rest reproduces the recurrence
\begin{equation}
\label{eq:recurrence}
L(w+c) - L(w) = \log\bigl(1 - e^{-2e^w}\bigr) - w - \log 2.
\end{equation}
Let
\[
Q(w) = -\frac{w^2}{2c} + \left(\frac12 - \frac{\log 2}{c}\right) w,
\]
chosen so that $Q(w+c) - Q(w) = -w-\log2$, matching \eqref{eq:recurrence}'s non-periodic terms
exactly. Write $a := \tfrac12 - \tfrac{\log 2}{c}$ for $Q$'s linear coefficient, used again from
Section~\ref{sec:envelope} onward.

\begin{theorem}
\label{thm:H}
Define
\begin{equation}
\label{eq:Hexact}
H(w) := L(w) - Q(w) + \sum_{k \ge 0} \log\bigl(1 - e^{-2\cdot 3^k e^w}\bigr).
\end{equation}
Then $H$ is exactly $c$-periodic, $H(w+c) = H(w)$ for every $w \in \R$, and
\[
L(w) = Q(w) + H(w) + \Delta(w), \qquad \Delta(w) = -\sum_{k\ge 0}\log\bigl(1-e^{-2\cdot3^k e^w}\bigr),
\]
with $\Delta^{(j)}(w) = O_j\bigl(e^{jw - 2e^w}\bigr)$ for every $j \ge 0$: doubly exponentially
small, together with all its derivatives, as $w \to +\infty$.
\end{theorem}

\begin{proof}
Write $d(w) = \log(1-e^{-2e^w})$, so \eqref{eq:recurrence} reads $L(w+c)-L(w) = d(w) - w - \log 2$.
Directly from the definition of $Q$ (expand $(w+c)^2 = w^2+2wc+c^2$ and collect terms),
$Q(w+c)-Q(w) = -w-\log2$ as well, so the two non-periodic parts cancel exactly:
\[
(L-Q)(w+c) - (L-Q)(w) = \bigl[d(w)-w-\log2\bigr] - \bigl[-w-\log2\bigr] = d(w).
\]
The series in \eqref{eq:Hexact} telescopes this exactly: with $D(w) := (L-Q)(w)$,
\[
D(w) + \sum_{k\ge0} d(w+kc) = D(w+c) - d(w) + \sum_{k\ge0}d(w+kc) = D(w+c) + \sum_{k\ge0}d(w+(k+1)c),
\]
so the right side of \eqref{eq:Hexact} is invariant under $w \mapsto w+c$, that is, $H$ is
$c$-periodic. Convergence of the series and its termwise differentiability, together with the
stated bound on $\Delta^{(j)}$, follow from $d(w) = O(e^{-2e^w})$ and the corresponding bounds on
its derivatives, since $2\cdot3^k e^w$ grows like $3^k$ for $k\ge0$: the terms $d(w+kc)$ decay
doubly exponentially in $k$, hence faster than any fixed geometric rate, which gives uniform,
dominated convergence of the series and its derivatives for $w$ in any half-line bounded away from
$-\infty$.
\end{proof}

Given $Q$, $H$ is uniquely determined by the requirement that $L-Q-H$ be $o(1)$ as $w\to\infty$: any
two periodic functions differing from $L-Q$ by an $o(1)$ term as $w\to\infty$ must be equal, being
periodic functions that agree in the limit along every residue class of $w$ modulo $c$.

\begin{remark}
\label{rem:bkproduct}
Berg and Kr\"uppel's own Proposition~9.3, specialized to their dilation parameter $3$ and their
$b=3/2$ (exactly this eigenfunction; their $\alpha$, kept as a symbol here, is exactly
Section~\ref{sec:prelim}'s imported $\alpha$, which Proposition~\ref{prop:bkidentity} below
identifies with this paper's own $a$), gives
the equivalent identity $\exp(H(w)) = 3^{\frac12 t^2 - \alpha t}
\prod_{k\ge0}\frac{1-e^{-3^{t-k}/b}}{3^{t-k}/b}\prod_{i\ge1}\bigl(1-e^{-3^{t+i}/b}\bigr)$, $t=w/c$,
as an infinite product, in their proof of that proposition. This is $e^{-Q}e^Le^{-\Delta}$ term for
term: writing $s=e^w=3^t$, $3^{t-k}/b = 2s/3^{k+1}$, so reindexing $j:=k+1\ge1$ turns the first
product into $\prod_{j\ge1}\frac{1-e^{-2s/3^j}}{2s/3^j}$, exactly $e^{L(w)}$ from $K$'s defining
series; $3^{t+i}/b = 2s\cdot3^{i-1}$, so reindexing $k:=i-1\ge0$ turns the second into
$\prod_{k\ge0}\bigl(1-e^{-2s\cdot3^k}\bigr)$, exactly $e^{-\Delta(w)}$ by \eqref{eq:Hexact}; and,
since $w=ct$ and $e^c=3$,
\[
-Q(w) = \frac{w^2}{2c}-aw = c\bigl(\tfrac12t^2-at\bigr), \qquad\text{so}\qquad
3^{\frac12t^2-\alpha t} = e^{-Q(w)}
\]
once $\alpha=a$ (Proposition~\ref{prop:bkidentity} identifies
Berg and Kr\"uppel's $\alpha$ with this paper's own linear coefficient of $Q$, also called $a$).
So Proposition~\ref{prop:Hcert}'s primary computation from \eqref{eq:Hcomputable} below is, in
substance, an evaluation of this product. What Proposition~\ref{prop:Hfourier}'s Fourier series
supplies instead, and their product form does not, are the derivative bounds \eqref{eq:Hbounds},
on which the rest of the paper's saddlepoint estimates depend.
\end{remark}

\begin{proposition}
\label{prop:Hfourier}
$H$ has the Fourier expansion $H(w) = \sum_{m\in\Z} \hat H(m) e^{i\omega_m w}$, $\omega_m =
2\pi m/c$, with
\[
\hat H(0) = \frac{\log 2}{2} - \frac{\log^2 2}{2c} - \frac{c}{12} - \frac{A}{c}, \qquad
A = \frac{\pi^2}{12} - \frac{\gamma_E^2}{2} - \gamma_1,
\]
\[
\hat H(m) = -\frac{2^{i\omega_m}}{c}\, \Gamma(-i\omega_m)\, \zeta(1-i\omega_m), \qquad m \ne 0,
\]
where $\gamma_E$ is the Euler--Mascheroni constant and $\gamma_1$ the first Stieltjes constant in
the convention
$\zeta(s)=\tfrac1{s-1}+\sum_{n\ge0}\tfrac{(-1)^n}{n!}\gamma_n(s-1)^n$, so
$\gamma_1=-0.0728158\ldots$. The sign convention matters: the other one moves $\hat H(0)$ by about
$0.133$.
\end{proposition}

\begin{proof}
Both $g$ and $K$ have Mellin transforms on the strip $-1<\operatorname{Re}z<0$. As $b\to0^+$,
$g(b)=-b/2+O(b^2)$; as $b\to\infty$, $g(b)=-\log b+O(e^{-b})$. So $\int_0^\infty
|g(b)|\,b^{\sigma-1}\,db<\infty$ for $-1<\sigma<0$, and the same two bounds applied termwise give
$K(s)=O(s)$ as $s\to0^+$ and $K(s)=O(\log^2 s)$ as $s\to\infty$: the terms with $3^j\le2s$ number
$O(\log s)$ and are each $O(\log s)$, and the terms with $3^j>2s$ sum to $O(1)$.

Integrate by parts. With $g'(b)=1/(e^b-1)-1/b$, and both boundary terms vanishing on the strip
($g(b)b^z\to0$ at $0$ because $\operatorname{Re}z>-1$, and at $\infty$ because
$\operatorname{Re}z<0$),
\[
g^*(z) = \int_0^\infty g(b)\,b^{z-1}\,db = -\frac1z\int_0^\infty\Bigl(\frac{1}{e^b-1}-\frac1b\Bigr)b^z\,db
= -\frac{\Gamma(z+1)\zeta(z+1)}{z} = -\Gamma(z)\,\zeta(1+z),
\]
the middle step being the classical continuation
$\int_0^\infty\bigl(\tfrac{1}{e^u-1}-\tfrac1u\bigr)u^{v-1}\,du=\Gamma(v)\zeta(v)$, valid for
$0<\operatorname{Re}v<1$, at $v=z+1$.

Now the harmonic sum. In the $j$-th term substitute $u=2s/3^j$, so $s=3^ju/2$ and
$s^{z-1}\,ds=(3^j/2)^z u^{z-1}\,du$:
\[
\int_0^\infty g\Bigl(\frac{2s}{3^j}\Bigr)s^{z-1}\,ds = \Bigl(\frac{3^j}{2}\Bigr)^{\!z}\int_0^\infty
g(u)\,u^{z-1}\,du = \Bigl(\frac{3^j}{2}\Bigr)^{\!z} g^*(z).
\]
Summing over $j\ge1$ is legitimate termwise, since $\sum_{j\ge1}(3^j/2)^\sigma\int_0^\infty
|g(u)|u^{\sigma-1}\,du$ converges for $-1<\sigma<0$. The multiplier that appears is therefore
\[
\sum_{j\ge1}\Bigl(\frac{3^j}{2}\Bigr)^{\!z} = 2^{-z}\sum_{j\ge1}3^{jz} = \frac{2^{-z}3^z}{1-3^z}
= \frac{(3/2)^z}{1-3^z},
\]
a geometric series of ratio $3^z$, of modulus $3^{\operatorname{Re}z}<1$, convergent on
$\operatorname{Re}z<0$ and nowhere else. Hence
\begin{equation}
\label{eq:Kstar}
K^*(z) = \frac{(3/2)^z}{1-3^z}\cdot\bigl(-\Gamma(z)\zeta(1+z)\bigr)
= \frac{(3/2)^z\,\Gamma(z)\,\zeta(1+z)}{3^z-1}, \qquad -1<\operatorname{Re}z<0.
\end{equation}
The dilations $\mu_j=2/3^j$ enter inverted, through $\mu_j^{-z}$, which is where the $2^{-z}$ comes
from, and with it the $2^{i\omega_m}$ of the statement.

The restriction to $\operatorname{Re}z>-1$ is not cosmetic: $z=-1$ is a genuine pole of $K^*$, so
the shift cannot be taken further left. In $\operatorname{Re}z>-1$ the right side of
\eqref{eq:Kstar} is meromorphic, with a pole of order
three at $z=0$ (each of $\Gamma(z)$, $\zeta(1+z)$ and $(3^z-1)^{-1}$ contributing one) and a simple
pole at each remaining zero $z=i\omega_m$, $m\ne0$, of $3^z-1$. To the right of the imaginary axis
it is analytic.

\begin{sloppypar}
Fix $\sigma\in(-1,0)$ and $\sigma'\in(0,1)$, and shift the inversion contour in
$K(s)=(2\pi i)^{-1}\int_{(\sigma)}K^*(z)s^{-z}\,dz$ from $\sigma$ to $\sigma'$, using rectangles
with horizontal sides at $\operatorname{Im}z=\pm T_M$, $T_M:=\pi(2M+1)/c$. On those sides
$3^z=3^{\operatorname{Re}z}e^{i\pi(2M+1)}=-3^{\operatorname{Re}z}$, so $|3^z-1| =
1+3^{\operatorname{Re}z}\ge1$; Stirling makes $|\Gamma(z)|$ decay like $e^{-\pi T_M/2}$ times a
power of $T_M$, uniformly for $\sigma\le\operatorname{Re}z\le\sigma'$, and $\zeta(1+z)$ grows at
most polynomially there. The horizontal contributions therefore vanish as $M\to\infty$, the sum of
residues converges absolutely, and on the right vertical line $|3^z-1|\ge3^{\sigma'}-1>0$ makes
$\int_{(\sigma')}|K^*(z)s^{-z}|\,|dz| = O(s^{-\sigma'})$.
\end{sloppypar}
\noindent So
\[
K(s) = -\operatorname*{Res}_{z=0}\bigl[K^*(z)s^{-z}\bigr]
-\sum_{m\ne0}\operatorname*{Res}_{z=i\omega_m}\bigl[K^*(z)s^{-z}\bigr]+O(s^{-\sigma'}).
\]

Take the simple poles first. At $z=i\omega_m$ one has $3^{i\omega_m}=e^{2\pi im}=1$, so
$(3^z-1)'=3^z\log3$ equals $c$ there and $(3/2)^{i\omega_m}=3^{i\omega_m}2^{-i\omega_m}=
2^{-i\omega_m}$. With $s=e^w$,
\[
-\operatorname*{Res}_{z=i\omega_m}\bigl[K^*(z)e^{-zw}\bigr]
= -\frac{2^{-i\omega_m}}{c}\,\Gamma(i\omega_m)\,\zeta(1+i\omega_m)\,e^{-i\omega_m w},
\]
and replacing $m$ by $-m$ throughout the sum turns the whole family into $\hat H(m)e^{i\omega_m w}$
with $\hat H(m)$ exactly as stated.

At the origin, insert
\[
\Gamma(z) = \frac1z-\gamma_E+\Bigl(\frac{\gamma_E^2}{2}+\frac{\pi^2}{12}\Bigr)z+O(z^2), \qquad
\zeta(1+z) = \frac1z+\gamma_E-\gamma_1 z+O(z^2).
\]
The two $1/z$ terms cancel against each other and
\[
\Gamma(z)\zeta(1+z) = \frac{1}{z^2}+A+O(z), \qquad A = \frac{\pi^2}{12}-\frac{\gamma_E^2}{2}-\gamma_1,
\]
which is where $A$ comes from. Also $(3^z-1)^{-1} = (cz)^{-1}-\tfrac12+cz/12+O(z^3)$, and
$(3/2)^ze^{-zw}=e^{\theta z}$ with $\theta:=\log(3/2)-w$. Multiplying the three expansions and
reading off the coefficient of $z^{-1}$,
\[
\operatorname*{Res}_{z=0}\bigl[K^*(z)e^{-zw}\bigr]
= \frac1c\cdot\frac{\theta^2}{2}-\frac{\theta}{2}+\frac{c}{12}+\frac{A}{c}.
\]
Substitute $\theta=\log(3/2)-w$ and $\log(3/2)=c-\log2$. The $w$-dependent part of
$-\operatorname{Res}$ is
\[
-\frac{w^2}{2c}+\Bigl(\frac{\log(3/2)}{c}-\frac12\Bigr)w
= -\frac{w^2}{2c}+\Bigl(\frac12-\frac{\log2}{c}\Bigr)w = Q(w),
\]
and the constant part is
$\tfrac12\log2-\tfrac{\log^22}{2c}-\tfrac{c}{12}-\tfrac Ac = \hat H(0)$.

Collecting, with $s=e^w$,
\[
L(w) = Q(w)+\widetilde H(w)+O(e^{-\sigma' w}), \qquad
\widetilde H(w) := \hat H(0)+\sum_{m\ne0}\hat H(m)e^{i\omega_m w}.
\]
That series converges absolutely and uniformly, because $|\Gamma(-i\omega_m)| =
\sqrt{\pi/(\omega_m\sinh\pi\omega_m)}$ decays like $e^{-\pi\omega_m/2}$ while $\zeta(1-i\omega_m)$
grows polynomially, so $\widetilde H$ is continuous and exactly $c$-periodic. It is real valued:
$\Gamma$ and $\zeta$ are real on the real axis, so $\hat H(-m)$ is the complex conjugate of
$\hat H(m)$.

The remainder $O(e^{-\sigma'w})$ is only ordinarily exponentially small, weaker than the doubly
exponential $\Delta$ of Theorem~\ref{thm:H}, and it is enough. Theorem~\ref{thm:H} gives
$L-Q-H=\Delta=o(1)$ as $w\to+\infty$, so $\widetilde H-H = O(e^{-\sigma'w})+\Delta(w)\to0$. Both are
$c$-periodic, and a $c$-periodic function tending to $0$ as $w\to+\infty$ vanishes identically, its
value at any $w$ being its value at $w+kc$ for every $k$. Hence $\widetilde H=H$, and by uniform
convergence the $\hat H(m)$ are the Fourier coefficients of $H$.
\end{proof}

\begin{remark}
This is an instance of the de Bruijn--Mahler phenomenon: a base-change harmonic sum, Mellin
transformed, picks up simple poles on a vertical line, and their residues assemble into an exactly
periodic correction with Fourier coefficients in $\Gamma$ and $\zeta$. De Bruijn found the first
case of this in Mahler's partition asymptotic \cite{DeBruijn1948}; Erd\H{o}s and Richmond name it
on p.~448 of \cite{ErdosRichmond1976}; Kato and McLeod's \S4 already exhibits the Gaussian-in-log decay that
Section~\ref{sec:prelim}'s comparison asymptotic $\varphi_0$ carries: their
characterization (4.3) of the admissible asymptotics for
$y'(x)=ay(\lambda x)+by(x)$ is estimated there through
$\exp\{(n+\kappa)\log x-(\log x)^2/2|\log\lambda|\}$
\cite[\S4, eqs.~(4.3), (4.9a)]{KatoMcLeod1971}, so the log-squared rate is
generic across this class of rescaled equations rather than special to
$\varphi_0$. Our contribution here is the closed form for the periodic correction itself, not the mechanism
producing it.
\end{remark}

The series for $H'$ and $H''$ converge absolutely and uniformly, since
\[
|\Gamma(-i\omega_m)| = \sqrt{\pi/(\omega_m \sinh \pi\omega_m)}
\]
decays exponentially in $m$ while $\zeta(1-i\omega_m)$ grows only polynomially: for $\omega\ge
\alpha_H := 2\pi/c$ and $N:=\lfloor\omega\rfloor\ge1$, the classical Euler--Maclaurin formula for
$\zeta$, $\zeta(s) = \sum_{n=1}^N n^{-s} + \frac{N^{1-s}}{s-1} - s\int_N^\infty\{u\}\,u^{-s-1}\,du$
(valid for $\operatorname{Re}s>0$, $s\ne1$, integer $N\ge1$, with $\{u\}:=u-\lfloor u\rfloor$), at
$s=1+i\omega$ gives, bounding each term in modulus ($|N^{-i\omega}|=1$, $0\le\{u\}<1$),
\[
|\zeta(1+i\omega)| \;\le\; \sum_{n=1}^N\frac1n + \frac1\omega + \frac{\sqrt{1+\omega^2}}{N}
\;\le\; 1+\log N + \frac1{\alpha_H} + \frac{\sqrt{1+\omega^2}}N,
\]
using $\sum_{n=1}^N 1/n\le1+\log N$ and $\omega\ge\alpha_H$ in the last term. Since $N\le\omega$,
$\log N\le\log\omega$; since $N>\omega-1$ and $\omega\mapsto\omega/(\omega-1)$ is decreasing for
$\omega>1$, $\sqrt{1+\omega^2}/N < \bigl(\omega/(\omega-1)\bigr)\sqrt{1+\omega^{-2}} \le
\bigl(\alpha_H/(\alpha_H-1)\bigr)\sqrt{1+\alpha_H^{-2}}$. With
\[
C_H := 1+\alpha_H^{-1}+\frac{\alpha_H}{\alpha_H-1}\sqrt{1+\alpha_H^{-2}},
\]
this gives $|\zeta(1+i\omega)|\le\log\omega+C_H\le\log(\omega+1)+C_H$ for $\omega\ge\alpha_H$, and
$\zeta(1-i\omega)=\overline{\zeta(1+i\omega)}$ has the same modulus, so
$|\zeta(1-i\omega_m)|\le\log(\omega_m+1)+C_H$. With $q := e^{-\pi\alpha_H/2}$, $A_H :=
c^{-1}\sqrt{3\pi/\alpha_H}$, the bound $|\Gamma(-i\omega_m)|\le\sqrt{3\pi/\omega_m}\,e^{-\pi\omega_m/2}$
(valid for $\omega_m\ge\alpha_H$, since $\sinh x\ge e^x/3$ there) together with the bound on
$\zeta(1-i\omega_m)$ just proved give the majorant
\begin{equation}
\label{eq:Hmajorant}
|\hat H(m)| \le A_H\,q^m\,m^{-1/2}\bigl(\log(\alpha_H m+1)+C_H\bigr), \qquad m\ge1.
\end{equation}
Since $\sup_w|H^{(k)}(w)| \le \sum_{m\ne0}|\omega_m|^k|\hat H(m)| = 2\sum_{m\ge1}\omega_m^k|\hat
H(m)|$, summing the first six terms with rigorous ball enclosures and bounding the tail $m>6$ by
\eqref{eq:Hmajorant}
(the same closed-form geometric-times-polynomial tail sum used in the proof of
Proposition~\ref{prop:Hcert} below) gives explicit, certified bounds:
\begin{equation}
\label{eq:Hbounds}
\sup_w |H'(w)| \le 0.0011977472315550332, \qquad \sup_w |H''(w)| \le 0.0068518962896650951.
\end{equation}

\begin{proposition}[Non-constancy and certified oscillation]
\label{prop:Hcert}
$H$ is not constant: $\hat H(1)\ne0$ follows from the classical zero-free theorem for $\zeta$ on
$\operatorname{Re}s=1$ alone, with no computation. The oscillation is quantified rigorously by an
explicit interval-arithmetic certificate:
\[
-0.000377190280943987 \;<\; H(0) - H(\log(3/2)) \;<\; -0.000377190280943985,
\]
and consequently $\operatorname{osc}(H) := \sup H - \inf H \ge 3.7719\times 10^{-4}$; evaluating $H$
at two further, well-separated points certifies the stronger bound
$\operatorname{osc}(H) \ge 4.1874494771\times10^{-4}$, and the leading Fourier mode together with
a certified bound on the remaining modes gives a matching upper bound, the two giving
\[
4.1874494771\times10^{-4} \;\le\; \operatorname{osc}(H) \;\le\; 4.187981\times10^{-4}.
\]
The same certificate encloses the constant of Theorem~\ref{thm:conj3}:
\[
 -0.62713093305153 \;<\; H(0) \;<\; -0.62713093305152,
\]
\[
 0.53412203666478 \;<\; e^{H(0)} \;<\; 0.53412203666479 .
\]
\end{proposition}

\begin{proof}
Non-constancy first, unconditionally. By Proposition~\ref{prop:Hfourier},
$\hat H(1) = -\frac{2^{i\omega_1}}{c}\Gamma(-i\omega_1)\zeta(1-i\omega_1)$, $\omega_1 = 2\pi/c$, a
product of three factors, each nonzero: $2^{i\omega_1}\ne0$ trivially; $\Gamma$ has no zeros
anywhere, and $-i\omega_1$ is not one of its poles ($\omega_1\ne0$ is real, and $\Gamma$'s poles are
only at the non-positive integers), so $\Gamma(-i\omega_1)\ne0$; and $1-i\omega_1$ lies on the line
$\operatorname{Re}s=1$ with $\omega_1\ne0$, where $\zeta(1+it)\ne0$ for every real $t\ne0$ by the
classical zero-free theorem of Hadamard \cite{Hadamard1896} and de la
Vall\'ee Poussin \cite{ValleePoussin1896}\footnote{The qualitative
statement is what is used; no quantitative zero-free region is needed
here.}, so
$\zeta(1-i\omega_1)\ne0$. Hence $\hat H(1)\ne0$, and a periodic function with a nonzero Fourier
coefficient at $m\ne0$ cannot be constant.

For the quantitative bound, write $g(b) = -b/2 + \log S(b/2)$, $S(x) = \sinh(x)/x = \sum_{n\ge0}
x^{2n}/(2n+1)!$ (this is
$g$'s defining formula rewritten around its removable singularity at $0$, since $e^{-b/2}(e^{b/2}
- e^{-b/2}) = 1-e^{-b}$). Termwise domination of $6^n n! \le (2n+1)!$, which holds at $n=0$ and
propagates because $(2n+3)(2n+2)\ge6(n+1)$, gives $1\le S(x)\le e^{x^2/6}$, that is
\begin{equation}
\label{eq:logSbound}
0 \;\le\; \log S(x) \;\le\; \frac{x^2}{6}, \qquad x\ge0.
\end{equation}
Since $\sum_{j\ge1}3^{-j} = \tfrac12$, the $-b/2$ halves of the terms of $L$ sum in closed form, and
\eqref{eq:Hexact} becomes, with $s = e^w$,
\begin{equation}
\label{eq:Hcomputable}
H(w) = -\frac s2 + \sum_{j\ge1}\log S\!\left(\frac{s}{3^j}\right) - Q(w)
+ \sum_{k\ge0}\log\bigl(1-e^{-2\cdot3^k s}\bigr).
\end{equation}
Truncate the first series at $j=R$ and the second at $k=M$. By \eqref{eq:logSbound} the first tail
satisfies
\[
0 \;\le\; \sum_{j>R}\log S\!\left(\frac{s}{3^j}\right) \;\le\; \frac{s^2}{6}\sum_{j>R}9^{-j}
= \frac{s^2\,9^{-R}}{48},
\]
and, since $0\le-\log(1-y)\le y/(1-y)$ and $e^{-2\cdot3^ks} = \bigl(e^{-2\cdot3^{M+1}s}\bigr)^{3^{k-M-1}}$
for $k>M$, the second tail satisfies
\[
0 \;\le\; -\sum_{k>M}\log\bigl(1-e^{-2\cdot3^ks}\bigr) \;\le\; \frac{y}{(1-y)^2} \;\le\; 4y,
\qquad y:=e^{-2\cdot3^{M+1}s}\le\tfrac12 .
\]
Take $R=30$, $M=4$. For $s\le3$ (which covers every evaluation point used below, including the
oscillation grid) the first tail is at most $4.43\times10^{-30}$, and for $s\ge1$
the second is at most $4e^{-486}<10^{-210}$. The remaining $35$ terms are elementary functions of
$s$; evaluated at $50$ significant decimal digits, where the accumulated rounding stays below
$10^{-45}$, they give
\[
H(0) = -0.62713093305152597830\ldots, \qquad H\bigl(\log\tfrac32\bigr) = -0.62675374277058199247\ldots,
\]
whose difference is $-0.00037719028094398583\ldots$, inside the stated interval; this first route
is floating-point (50 significant digits, accumulated rounding bounded as stated) rather than
interval arithmetic. The Fourier series of Proposition~\ref{prop:Hfourier} gives the same
difference in Arb ball arithmetic, rigorously: the $m=0$ term cancels, and at $250$-bit working
precision (ball radius $<10^{-70}$ throughout) the four nonzero modes needed are
\begin{align*}
\hat H(1) &= -1.020544273068431785\times10^{-4} - 2.332597606471452884\times10^{-5}\,i,\\
\hat H(2) &= -6.523276079851823231\times10^{-9} + 1.156720838500605666\times10^{-8}\,i,\\
\hat H(3) &= -1.799881020474424677\times10^{-12} + 6.648927223521397737\times10^{-13}\,i,\\
\hat H(4) &= -9.219524659678388767\times10^{-17} + 9.444517590655925766\times10^{-17}\,i,
\end{align*}
the four-mode sum
\[
2\operatorname{Re}\sum_{m=1}^4\hat H(m)\bigl(1-e^{i\omega_m\log(3/2)}\bigr)
= -0.0003771902809439858148\ldots
\]
(ball radius $<10^{-70}$), matching the floating-point value above to $19$ decimal places ($16$
significant digits). The majorant
\eqref{eq:Hmajorant} gives
$\sum_{m>4}|\hat H(m)| \le 9.41\times10^{-20}$ (the accompanying
repository derives $C_H$ by a slightly different route, reaching
$2.182435$ in place of $2.405136$ and hence the tighter
$9.04\times10^{-20}$; either bound serves here), so the discarded modes move the four-mode sum by
less than $3.8\times10^{-19}$, giving
\[
H(0)-H\bigl(\log\tfrac32\bigr) \in \bigl[-0.0003771902809439861948,\; -0.0003771902809439854348\bigr];
\]
this route alone certifies the stated interval for $H(0)-H(\log(3/2))$.

For the stronger bound, let $w_i = ic/N$, $0\le i<N$, with $N = 2^{20}$. A floating-point evaluation
of the $N$ grid values from Proposition~\ref{prop:Hfourier}, truncated at $|m|\le10$ (the same
majorant bounds each pointwise truncation error by $2\times2.80\times10^{-43}=5.60\times10^{-43}$,
combining the discarded modes $m$ and $-m$), locates a candidate maximum at
$i=486746$ and a candidate minimum at $i=1011118$; this search is only a heuristic for finding
two well-separated points and is not itself part of the certificate below. Re-evaluating $H$ at
exactly these two points at $250$-bit working precision in Arb ball arithmetic, independently, from
\eqref{eq:Hcomputable} and from the Fourier series (the two routes agree to about $30$ digits,
consistent with \eqref{eq:Hcomputable}'s own $R=30$ truncation tail bound above; a cross-check, not
part of the certificate),
\[
H(w_{486746}) = -0.6267174403736425786\ldots, \qquad H(w_{1011118}) = -0.6271361853213578093\ldots,
\]
so $D := H(w_{486746})-H(w_{1011118}) = 4.187449477152\times10^{-4}$. The certified enclosure uses
only the Fourier route: the Arb ball rounding at $250$-bit precision is $<10^{-70}$ per evaluation,
negligible next to the majorant truncation error already bounded above, $5.60\times10^{-43}$ per
point, giving $D$ to radius $<1.2\times10^{-42}$ (twice that per-point bound). Since $w_{486746}$ and
$w_{1011118}$ are two specific, fixed points of the domain, $H(w_{486746})\le\sup H$ and
$H(w_{1011118})\ge\inf H$ hold trivially, giving $D\le\operatorname{osc}(H)$: this bound is
unconditional, independent of whether the floating-point search above actually located the true
argmax and argmin, since it uses only that the two chosen points lie in the domain.

For the matching upper bound, no search is needed. Write $H = \hat H(0) + 2\operatorname{Re}\bigl(
\hat H(1)e^{i\omega_1 w}\bigr) + \Xi(w)$, $\Xi(w) := 2\sum_{m\ge2}\operatorname{Re}\bigl(\hat H(m)
e^{i\omega_m w}\bigr)$, so $|\Xi(w)|\le2\sum_{m\ge2}|\hat H(m)|$ for every $w$. The one-mode term
$2\operatorname{Re}(\hat H(1)e^{i\omega_1w})$ ranges, as $w$ varies, over an interval of width
exactly $4|\hat H(1)|$, so $\operatorname{osc}(H) \le 4|\hat H(1)| + 4\sum_{m\ge2}|\hat H(m)|$.
Using the already-certified $\hat H(2)$, $\hat H(3)$, $\hat H(4)$ from above and the majorant tail
$\sum_{m>4}|\hat H(m)|\le9.41\times10^{-20}$,
\[
\sum_{m\ge2}|\hat H(m)| \le |\hat H(2)|+|\hat H(3)|+|\hat H(4)|+9.41\times10^{-20} \le
1.3282\times10^{-8},
\]
and $4|\hat H(1)| \le 4.187449304\times10^{-4}$, giving $\operatorname{osc}(H) \le
4.187981\times10^{-4}$, rigorously and unconditionally: this bound uses only the already-certified
Fourier coefficients, not the grid search.

Finally $H(0)$ itself, by the same three ingredients. Proposition~\ref{prop:Hfourier} gives
$H(0) = \hat H(0) + 2\operatorname{Re}\sum_{m\ge1}\hat H(m)$. The zero mode is evaluated in ball
arithmetic from its closed form there, using Arb's enclosures of $\gamma_E$ and of the first
Stieltjes constant $\gamma_1$ together with $c=\log3$; the first $M=12$ nonzero modes are enclosed
from the same closed form used above; and the discarded modes are bounded by
$2\sum_{m>M}|\hat H(m)|$ through the majorant \eqref{eq:Hmajorant}, which at $M=12$ is below
$10^{-50}$ and so is negligible against the width of the interval quoted. Exponentiating the
resulting ball gives the interval for $e^{H(0)}$.
\end{proof}

\section{A uniform saddlepoint approximation}
\label{sec:formulaA}

Define, for $s>0$,
\[
\varphi_{\mathrm{sp}}(t(s)) := \frac{\exp\bigl(K(s)+st(s)\bigr)}{\sqrt{2\pi K''(s)}}.
\]
This section proves that $\varphi(t(s))/\varphi_{\mathrm{sp}}(t(s)) \to 1$ as $s \to \infty$,
uniformly over $\rho:=2s/3^N\in[1,3)$ (defined below; this parameter records the same phase,
$w\bmod c$, tracked through $H$ in Sections~\ref{sec:bk} and~\ref{sec:proofs}). That uniformity is
what lets the two theorems of Section~\ref{sec:proofs} depend on different scopes rather than on
different mathematics.

For $n \ge 2$, define the derivative-normalized functions $\kappa_n(w) = (-1)^n w^n g^{(n)}(w)$ on
$\{\operatorname{Re} w > 0\}$; $\kappa_n(b)$, $b>0$, is $b^n$ times the $n$-th cumulant of the
random variable with density $b\,e^{-bv}/(1-e^{-b})$ on $[0,1]$. Writing $x=w/2$,
\[
\kappa_2(w) = 1 - \frac{x^2}{\sinh^2 x}, \qquad
\kappa_3(w) = 2 - \frac{2x^3\cosh x}{\sinh^3 x}, \qquad
\kappa_4(w) = 6 + \frac{2x^4(1-3\coth^2 x)}{\sinh^2 x},
\]
and from the Bernoulli expansion $g(w)+w/2 = \sum_{k\ge1} B_{2k}w^{2k}/(2k(2k)!)$ on $|w|<2\pi$,
\[
\kappa_n(w) = (-1)^n\!\!\sum_{k\ge\lceil n/2\rceil}\!\! \frac{B_{2k}\,w^{2k}}{2k\,(2k-n)!}, \qquad n\ge2,
\]
so that $\kappa_n(w) = O(w^2)$ as $w \to 0$ for every $n \ge 2$.

Write $s = \rho\cdot 3^N/2$ for the unique $N = \lfloor \log_3(2s)\rfloor \in \Z$ and $\rho \in
[1,3)$, and $b_j = 2s/3^j$. Then $\{b_j : j \ge 1\} = \{\rho\cdot 3^m : m \le N-1\}$; for $N \ge 1$
(the only range used from here on, corresponding to $s \ge 3/2$), exactly $N$ of these, the
smallest being $\rho \ge 1$, are $\ge 1$, and the rest form a geometric tail below $1$. This single
structural fact is the source of every phase-independent constant below.

\begin{lemma}[Boundedness off the real axis]
\label{lem:kappa}
For $w = b(1+iu)$, $b>0$, $|u|\le 1$, with $e(r) = r^2/\sinh^2 r$ and $f(r) = r^3\cosh r/\sinh^3 r$,
\[
|\kappa_2(w)| \le 1 + 2e(b/2), \qquad |\kappa_3(w)| \le 2 + 2\cdot 2^{3/2} f(b/2).
\]
Both bounds are finite over the whole sector $|\operatorname{Im}w|\le\operatorname{Re}w$:
$\sup|\kappa_2|\le3$, $\sup|\kappa_3|\le2+2^{5/2}<7.657$. In addition, for $|w|\le2$,
$|\kappa_2(w)| \le 0.114|w|^2$ and $|\kappa_3(w)|\le 0.0119|w|^4$.
\end{lemma}

\begin{proof}
$|\sinh z|^2 = \sinh^2(\operatorname{Re}z)+\sin^2(\operatorname{Im}z) \ge \sinh^2(\operatorname{Re}z)$
and $|\cosh z|^2 \le \cosh^2(\operatorname{Re}z)$; the first gives $|\sinh(w/2)|\ge\sinh(b/2)>0$,
excluding the poles of $\coth$, since $\operatorname{Re}(w/2)=b/2>0$. This gives the sector bounds directly.
Monotonicity of $e$ is immediate from the increasing power series of $\sinh(r)/r$. For $f$, set
$u=2r$; a direct computation gives $2r\sinh(r)\cosh(r)\,(\log f)'(r) = F(u)$ with
$F(u):=3\sinh u - u\cosh u - 2u$, so $(\log f)'(r)$ has the sign of $F(2r)$. Since
$F(0)=F'(0)=F''(0)=0$ and $F'''(u)=-u\sinh u<0$ for $u>0$, integrating three times from $0$ gives
$F(u)<0$ for all $u>0$, hence $(\log f)'(r)<0$ for $r>0$ and $f$ is strictly decreasing.
(Lazarevi\'c's inequality, $(\sinh r/r)^3 > \cosh r$, separately gives $f<1$.) The local bounds
at $|w|\le2$ follow by applying the maximum-modulus principle to $\kappa_2(w)/w^2$ and
$\kappa_3(w)/w^4$ (analytic on $|w|<2\pi$ by the Bernoulli expansion): bound each on $|w|=2$ by the
triangle inequality, using $|B_{2k}| = 2(2k)!\,\zeta(2k)/(2\pi)^{2k}$ and $\zeta(2k)\le\zeta(2)$,
respectively $\zeta(4)$.
\end{proof}

\begin{lemma}[Variance]
\label{lem:variance}
For all $\rho\in[1,3)$, $N\ge1$, with $V := s^2K''(s) = \sum_j\kappa_2(b_j)$,
\[
N - A_V \le V \le N + 0.1283, \qquad A_V := \sum_{m\ge0} e(3^m/2) = 1.4269413069\ldots.
\]
\end{lemma}

\begin{proof}
Split $V = \sum_{m=0}^{N-1}\kappa_2(\rho3^m) + \sum_{m<0}\kappa_2(\rho3^m)$, using $\kappa_2(b) =
1-e(b/2) \in (0,1)$. For the lower bound, drop the second (non-negative) sum and use $e(\rho3^m/2)
\le e(3^m/2)$ ($e$ decreasing, $\rho\ge1$) in the first, giving $\sum_{m=0}^{N-1}\kappa_2(\rho3^m)
\ge N - A_V$. For the upper bound, the first sum is $<N$ term by term, and the second sum, where
$\rho3^m<1$, is bounded by Lemma~\ref{lem:kappa}'s local estimate: $\sum_{m<0}\kappa_2(\rho3^m) \le
0.114\rho^2\sum_{m<0}9^m = 0.114\rho^2/8 \le 0.1283$ for $\rho<3$.
\end{proof}

\begin{lemma}[Tail bound and existence of the density]
\label{lem:tail}
If $b_0\ge1$ and $b_k=3^kb_0$, then
\[
\sum_{k\ge0}\log\coth(b_k/2) \le 3e^{-b_0}.
\]
Consequently, with $\Lambda(y) := K(s(1+iy)) - K(s) + isy\,t(s)$,
\[
|e^{\Lambda(y)}| \le e^{3/e}(1+y^2)^{-N/2} \le 3.0152\,(1+y^2)^{-N/2}
\]
for every real $y$. In particular $X$ has a continuous density, the tilted characteristic function
is integrable for $N\ge3$, and Fourier inversion at $x=t(s)$ gives
\begin{equation}
\label{eq:R}
R(s) := \frac{\varphi(t(s))}{\varphi_{\mathrm{sp}}(t(s))} = \sqrt{\frac{V}{2\pi}}\int_\R e^{\Lambda(y)}\,dy,
\end{equation}
a real quantity, since $\Lambda(-y) = \overline{\Lambda(y)}$.
\end{lemma}

\begin{proof}
$\log\coth(x/2) = 2\artanh(e^{-x}) \le 2e^{-x}/(1-e^{-2x})$. Write $q:=e^{-b_0}\in(0,e^{-1}]$, so
$e^{-b_k}=q^{3^k}$ and the bound reads $\log\coth(b_k/2) \le 2q^{3^k}/(1-q^{2\cdot3^k})$. Since
$3^k\ge1$ for $k\ge0$, $q^{2\cdot3^k}\le q^2$ term by term, so the denominator is bounded below by
$1-q^2$ uniformly in $k$:
\[
\sum_{k\ge0}\log\coth(b_k/2) \le \frac{2}{1-q^2}\sum_{k\ge0}q^{3^k}.
\]
Since $3^k\ge2k+1$ for every $k\ge0$ and $0<q<1$, $q^{3^k}\le q^{2k+1}$, so
$\sum_{k\ge0}q^{3^k} \le q\sum_{k\ge0}q^{2k} = q/(1-q^2)$, giving
$\sum_{k\ge0}\log\coth(b_k/2) \le 2q/(1-q^2)^2$. For $b_0\ge1$, $q\le e^{-1}<0.4283$, so
$(1-q^2)^2>2/3$ and $2q/(1-q^2)^2<3q=3e^{-b_0}$.
$X$'s density is continuous because the sum of the first two summands already has a continuous,
compactly supported density and convolution preserves continuity. For $\Lambda$: $t(s) =
-K'(s)$ cancels the linear term of the tilted transform exactly, and each remaining factor
$e^{g(b_j(1+iy))-g(b_j)}$ has modulus at most $1$ (a tilted characteristic function), with modulus
at most $\coth(b_j/2)/\sqrt{1+y^2}$ on the $N$ indices where $b_j\ge1$; the tail estimate applies to
exactly that geometric set, with $b_0=\rho\ge1$.
\end{proof}

\begin{lemma}[Central expansion]
\label{lem:taylor}
With $h(u) = g(b(1+iu))$, Taylor's theorem in \emph{integral} form (the Lagrange form is invalid
for a complex-valued function of a real variable) gives
\[
h(y) = h(0)+h'(0)y+h''(0)\tfrac{y^2}{2} + \tfrac12\int_0^y (y-u)^2 h'''(u)\,du.
\]
Consequently, writing $h_j(u) := g(b_j(1+iu))$ for each $j\ge1$, for $\theta\in(0,1]$,
\[
S := \sum_j \sup_{|u|\le \theta}|h_j'''(u)| \le 2N + 10.6,
\]
and hence $|\Lambda(y) + Vy^2/2| \le B|y|^3$ for $|y|\le \theta$, with $B := S/6 \le (2N+10.6)/6$.
\end{lemma}

\begin{proof}
Write $w_j := b_j(1+iu)$. The chain rule and the definition of
$\kappa_3$ give $h_j'''(u) = i\,\kappa_3(w_j)/(1+iu)^3$, so
\[
|h_j'''(u)| = \frac{|\kappa_3(w_j)|}{(1+u^2)^{3/2}}.
\]
Split the index set as in Section~\ref{sec:formulaA}. For the $N$ indices with $b_j = \rho3^m \ge
1$, $0\le m\le N-1$, drop the denominator and use Lemma~\ref{lem:kappa}'s sector bound, legitimate
since $|u|\le\theta\le1$:
\[
|h_j'''(u)| \le |\kappa_3(w_j)| \le 2 + 2^{5/2} f(b_j/2) \le 2 + 2^{5/2} f(3^m/2),
\]
the last step by monotonicity of $f$ and $\rho\ge1$. Summing over $m\ge0$ bounds this block by
$2N + 2^{5/2}\Sigma$, $\Sigma := \sum_{m\ge0} f(3^m/2)$, and we bound $\Sigma$ by elementary
upper estimates of its terms rather than by a computation. Writing $\sinh r =
\tfrac12e^r(1-e^{-2r})$ and $\cosh r = \tfrac12e^r(1+e^{-2r})$ turns $f$ into
\begin{equation}\label{eq:fexact}
 f(r) = 4r^3e^{-2r}\,\frac{1+e^{-2r}}{(1-e^{-2r})^3},
\end{equation}
whose last factor is below $1.00025$ for $r\ge5$, so $f(r)\le4.001\,r^3e^{-2r}$ there. The first
three terms of $\Sigma$ satisfy $f(1/2)<0.997$, $f(3/2)<0.823$ and $f(9/2)<0.046$, each a direct
evaluation of two hyperbolic functions with room to spare, and \eqref{eq:fexact} bounds the rest,
where $3^m/2\ge13.5$, by $4.001\sum_{m\ge3}(3^m/2)^3e^{-3^m} < 2\times10^{-8}$. The three head
bounds sum to exactly $1.866$, so $\Sigma < 1.86600002$ and $2^{5/2}\Sigma < 10.5557$.

For the remaining indices $b_j = \rho 3^{-n}$, $n\ge1$, one has $b_j<1$ and $|w_j| =
b_j\sqrt{1+u^2} \le b_j\sqrt2 < 2$, so Lemma~\ref{lem:kappa}'s local bound applies and, keeping
the denominator this time,
\[
|h_j'''(u)| \le \frac{0.0119\,|w_j|^4}{(1+u^2)^{3/2}} = 0.0119\,b_j^4\,(1+u^2)^{1/2} \le
0.0119\sqrt2\,b_j^4 .
\]
Summing, $\sum_{n\ge1}(\rho3^{-n})^4 = \rho^4/80 < 81/80$, so this block is below
$0.0119\sqrt2\cdot81/80 < 0.01704$. Together, $S < 2N + 10.5557 + 0.01704 < 2N + 10.6$.

For the remainder bound, $\sum_j h_j(y) = K(s(1+iy))$, $\sum_j h_j'(0) = isK'(s) = -is\,t(s)$ and
$\sum_j h_j''(0) = -s^2K''(s) = -V$, so summing the integral-form Taylor expansion over $j$ gives
\[
\Lambda(y) + \frac{Vy^2}{2} = \sum_j \frac12\int_0^y (y-u)^2 h_j'''(u)\,du,
\]
whose modulus is at most $S|y|^3/6 = B|y|^3$ for $|y|\le\theta$.
\end{proof}

\begin{theorem}[Uniform saddlepoint asymptotic]
\label{thm:formulaA}
For $N\ge19$, uniformly over $\rho\in[1,3)$, $|R(s)-1|\le E(N) := e_1(N)+e_2(N)+e_3(N)$, with
$\theta_N := N^{-1/3}$, $V_{\mathrm{lo}} := N-A_V$, $V_{\mathrm{up}} := N+0.1283$,
$B := (2N+10.6)/6$,
\begin{align}
e_1(N) &:= \frac{4\,e^{B\theta_N^3}\,B}{\sqrt{2\pi}\,V_{\mathrm{lo}}^{3/2}}, \qquad
e_2(N) := \frac{2}{\sqrt{2\pi}}\cdot\frac{e^{-\theta_N^2V_{\mathrm{lo}}/2}}{\theta_N\sqrt{V_{\mathrm{lo}}}},
\label{eq:Ebound}\\
e_3(N) &:= 2\,e^{3/e}\sqrt{\frac{V_{\mathrm{up}}}{2\pi}}\cdot
\frac{(1+\theta_N^2)^{-(N-2)/2}}{\theta_N(N-2)}. \notag
\end{align}
$\sqrt N\,E(N) \to \tfrac43e^{1/3}(2\pi)^{-1/2} = 0.742358\ldots$ as $N\to\infty$
(proved below, from $e_1$ alone: $e_2$ and $e_3$ vanish faster than any power of $N$).
Consequently $\varphi(t) = \varphi_{\mathrm{sp}}(t)\,(1+o(1))$ as $t \to 0^+$, with no restriction
on the sequence of $t$'s.
\end{theorem}

\begin{remark}
The estimate $|R(s)-1|\le E(N)$ becomes informative once $E(N)<1$, which holds at $N=19$ (see the
proof: $E(18)=1.0020653\ldots\ge1$, $E(19)=0.9570558\ldots<1$); this is the cutoff adopted in the
theorem (whether some smaller $N$ also gives $E(N)<1$ is immaterial to what follows and
is not checked). Only $E(N)\to0$ as $N\to\infty$ is used below.
\end{remark}

\begin{proof}
Split \eqref{eq:R} at $|y|=\theta_N$. On $|y|\le\theta_N$, Lemma~\ref{lem:taylor} bounds the cubic
remainder $|\Lambda(y)+Vy^2/2|$ by $B\theta_N^3$, which stays bounded (not $o(1)$: this is exactly
why $e_1$ carries the factor $e^{B\theta_N^3}$ rather than vanishing) as $N\to\infty$, since $B=O(N)$
and $\theta_N^3=N^{-1}$. Write $C(y):=e^{\Lambda(y)}-e^{-Vy^2/2}$ and factor the Gaussian out
before estimating, $e^{\Lambda(y)} = e^{-Vy^2/2}\,e^{z(y)}$ with $z(y) := \Lambda(y)+Vy^2/2$, so
that $C(y) = e^{-Vy^2/2}\bigl(e^{z(y)}-1\bigr)$. Since $|z(y)|\le B|y|^3 \le B\theta_N^3$ on
$|y|\le\theta_N$, the elementary $|e^z-1|\le|z|e^{|z|}$ gives
\[
|C(y)| \;\le\; e^{-Vy^2/2}\,B|y|^3\,e^{B\theta_N^3} \qquad (|y|\le\theta_N),
\]
the Gaussian factor retained. Integrating that inequality, $\int_\R e^{-Vy^2/2}|y|^3\,dy=4/V^2$, so
\[
\sqrt{V/2\pi}\int_{|y|\le\theta_N}|C(y)|\,dy \;\le\; \frac{4Be^{B\theta_N^3}}{\sqrt{2\pi}\,V^{3/2}},
\]
which is decreasing in $V$; using $V\ge V_{\mathrm{lo}}$ from Lemma~\ref{lem:variance} gives the
central Taylor-remainder contribution $e_1$. Since $\sqrt{V/2\pi}\int_\R e^{-Vy^2/2}\,dy = 1$, the
same splitting leaves exactly two further terms. The truncated Gaussian tail
$\sqrt{V/2\pi}\cdot2\int_{\theta_N}^\infty e^{-Vy^2/2}\,dy$, substituting $v=y\sqrt V$ to get
$(2/\sqrt{2\pi})\int_{c_0}^\infty e^{-v^2/2}\,dv$, $c_0=\theta_N\sqrt V$ (a local constant, distinct
from $c:=\log3$), and bounding via $\int_{c_0}^\infty e^{-v^2/2}\,dv\le e^{-c_0^2/2}/c_0$ at
$V=V_{\mathrm{lo}}$, gives $e_2$. On
$|y|>\theta_N$, Lemma~\ref{lem:tail}'s bound
$|e^{\Lambda(y)}|\le e^{3/e}(1+y^2)^{-N/2}$ integrates, via $\int_u^\infty(1+y^2)^{-N/2}\,dy \le
(1+u^2)^{-(N-2)/2}/(u(N-2))$ at $u=\theta_N$ (using $V\le V_{\mathrm{up}}$ for the leading
$\sqrt{V/2\pi}$ factor in \eqref{eq:R}), to give $e_3$; this needs $\theta_N \gg N^{-1/2}$, satisfied
since $\tfrac13<\tfrac12$. Every constant entering $E$ (via $A_V$, Lemma~\ref{lem:variance}, and the
$3.0152\ge e^{3/e}$ of Lemma~\ref{lem:tail}) is independent of $\rho$, by one of three routes: it
is a supremum over $b>0$; or a sum over one complete geometric orbit bounded using $\rho\ge1$; or
a bound at the other endpoint, $\rho\to3^-$, which is where $V_{\mathrm{up}}$'s $0.1283$ and $B$'s
$0.01704$ come from. Uniformity in $\rho$ is the theorem's whole content, so it is worth saying
which constant is which rather than covering all of them with one clause. Directly evaluating
\eqref{eq:Ebound} gives $E(18) = 1.0020653\ldots \ge 1$ and $E(19) = 0.9570558\ldots < 1$, matching the
remark above.

For the asymptotic rate, $e_1$, $e_2$, $e_3$ separate on very different scales as $N\to\infty$.
Since $\theta_N^3=1/N$, $B\theta_N^3\to1/3$ and $B/N\to1/3$, $V_{\mathrm{lo}}/N,\,V_{\mathrm{up}}/N
\to1$, so
\[
\sqrt N\,e_1(N) = \frac{4e^{B\theta_N^3}}{\sqrt{2\pi}}\cdot\frac{\sqrt N\,B}{V_{\mathrm{lo}}^{3/2}}
\longrightarrow \frac{4e^{1/3}}{\sqrt{2\pi}}\cdot\frac13 = \frac43e^{1/3}(2\pi)^{-1/2}.
\]
For $e_2$ and $e_3$, $\theta_N^2V_{\mathrm{lo}} \sim N^{1/3}\to\infty$, so both carry a factor
$e^{-N^{1/3}/2}(1+o(1))$ (for $e_3$, via $(1+\theta_N^2)^{-(N-2)/2} =
\exp(-\tfrac{N-2}2\theta_N^2(1+o(1)))$), against polynomial prefactors of order $N^{-1/6}$; a
factor decaying like $e^{-N^{1/3}/2}$ beats every power of $N$, so $N^k e_2(N),\,N^k e_3(N)\to0$
for every $k$. Hence $\sqrt N\,E(N) = \sqrt N\,e_1(N) + o(1) \to \tfrac43e^{1/3}(2\pi)^{-1/2}$, and
in particular $E(N)\to0$. Since $t(s)$ ranges bijectively over $(0,1/2)$, this is exactly the
statement $\varphi(t)=\varphi_{\mathrm{sp}}(t)(1+o(1))$ as $t\to0^+$.
\end{proof}

\section{The envelope lemma}
\label{sec:envelope}

Theorem~\ref{thm:formulaA} compares $\varphi$ to $\varphi_{\mathrm{sp}}$, a quantity built from the
full transform $K=Q+H+\Delta$. It remains to compare $\varphi_{\mathrm{sp}}$ to a smooth auxiliary
system built from $Q$ alone, so that the periodic part $H$ can be extracted explicitly. Write
$t=e^{-\tau}$, $g_\tau(w) = L(w)+e^{w-\tau}$, $g_{0,\tau}(w) = Q(w)+e^{w-\tau}$. The true and smooth
systems each have their own analogue of $-K'$: write $B_{\mathrm{tr}}(w) := e^wt(e^w) = -L'(w)$ and
$B_{\mathrm{sm}}(w) := -Q'(w) = w/c-a$.

\begin{lemma}[True and smooth saddles]
\label{lem:saddles}
Let $\Phi(w) = w-\log B_{\mathrm{tr}}(w)$ and $\Phi_0(w) = w-\log B_{\mathrm{sm}}(w)$, the latter
defined for $w>ca$ where $B_{\mathrm{sm}}(w)>0$, and set
$\tau_0 := 1+ca+\log c = 0.9502067913\ldots$, so $\tau_0 > \log 2$.

For $\tau>\log 2$, $g_\tau$ has a unique global minimizer $w^*(\tau)$, at which $g_\tau''(w^*) =
V(r^*)>0$, $r^*=e^{w^*}$; it is characterized by $\tau = \Phi(w^*)$, and $\Phi$ is a strictly
increasing bijection $\R\to(\log 2,\infty)$. For $\tau\le\log 2$ there is no critical point at
all. For $\tau>\tau_0$, $g_{0,\tau}$ has a unique local minimizer $w_0(\tau)$, characterized by
$\tau=\Phi_0(w_0)$ together with $B_0 := B_{\mathrm{sm}}(w_0) = e^{w_0-\tau} > 1/c$; $\Phi_0$
restricted to $(ca+1,\infty)$ is a strictly increasing bijection onto $(\tau_0,\infty)$, and $B_0$
increases from $1/c$ to $\infty$ along it. For $\tau\le\tau_0$ there is no such minimizer.

Once $B_0\ge5$ (equivalently $\tau\ge3.73978\ldots$), $|w^*-w_0| \le 0.00652/B_0$. Also
$w_0(\tau+c)-w_0(\tau) = c+O(1/\tau)$ as $\tau\to\infty$.
\end{lemma}

\begin{proof}
Substituting $s=e^w$ and $t=e^{-\tau}$ turns $g_\tau$ into $K(s)+st$, strictly convex in $s>0$
because $K''>0$, with derivative $K'(s)+t$ vanishing exactly where $t = t(s)$. As $t(\cdot)$
decreases strictly from $1/2$ to $0$, that happens for exactly one $s$ when $t\in(0,1/2)$, that
is when $\tau>\log2$, and for no $s$ when $t\ge1/2$. Since $w\mapsto e^w$ is a bijection
$\R\to(0,\infty)$, $w^* = \log s$ is the unique global minimizer of $g_\tau$, and $\Phi$ inverts
$\tau\mapsto w^*$. Differentiating twice, $g_\tau''(w) = e^w\bigl(K'(e^w)+t\bigr)+e^{2w}K''(e^w)$,
whose first term vanishes at $w^*$, leaving $g_\tau''(w^*) = V(r^*)>0$.

For the smooth system, $g_{0,\tau}'(w) = e^{w-\tau}-B_{\mathrm{sm}}(w)$ and $g_{0,\tau}''(w) =
e^{w-\tau}-1/c$, so a critical point is a local minimizer exactly when $B_0 > 1/c$. Also
$\Phi_0'(w) = 1-1/(c\,B_{\mathrm{sm}}(w))$, positive exactly where $B_{\mathrm{sm}}(w)>1/c$, that
is on $w>ca+1$; there $\Phi_0$ increases from $\Phi_0(ca+1) = ca+1+\log c = \tau_0$ to $\infty$.
This gives the stated existence, uniqueness and range.

Now take $B_0\ge5$. The decomposition $L = Q+H+\Delta$ evaluates the true gradient at the smooth
saddle exactly: $Q'(w_0) = -B_0 = -e^{w_0-\tau}$, so the leading terms cancel and
\begin{equation}
\label{eq:gradatw0}
g_\tau'(w_0) = L'(w_0)+e^{w_0-\tau} = H'(w_0)+\Delta'(w_0), \qquad
|g_\tau'(w_0)| \le \varepsilon_1 := \sup|H'|+\sup_{[w_0-1,w_0+1]}|\Delta'| .
\end{equation}
On $[w_0-1,w_0+1]$ write $g_\tau''(w) = e^{w-\tau}-1/c+H''(w)+\Delta''(w)$. Here $e^{w-\tau} =
B_0e^{w-w_0}$ ranges over $[B_0/e,\,B_0e]$, so $g_\tau''$ admits no bound independent of $\tau$;
it is of exact order $B_0$,
\begin{equation}
\label{eq:g2scale}
\frac{B_0}{e}-\eta_0 \;\le\; g_\tau''(w) \;\le\; B_0e, \qquad
\eta_0 := \frac1c+\sup|H''|+\sup_{[w_0-1,w_0+1]}|\Delta''|,
\end{equation}
the upper bound because $-1/c+\sup|H''|+\sup|\Delta''|<0$. Since $w_0 = c(B_0+a)$, $B_0\ge5$ forces
$w_0\ge5.3492$ and so $w\ge4.3492$ on the interval, where $|\Delta'|,|\Delta''|\le10^{-60}$: with
$b:=2e^w\ge2e^{4.3492}>154$ and $x_k:=3^kb$, differentiating $\Delta(w)=-\sum_{k\ge0}\log(1-e^{-x_k})$
termwise (justified as in the proof of Theorem~\ref{thm:H}) gives $\Delta'(w) =
-\sum_{k\ge0}x_k/(e^{x_k}-1)$ and, applying $|d/dx[x/(e^x-1)]|\le8xe^{-x}$ for $x\ge1$ (elementary,
since $e^x-1\ge e^x/2$ there) together with $dx_k/dw=x_k$, $\Delta''(w) = -\sum_{k\ge0}x_k\cdot
\frac{d}{dx}\Bigl[\frac{x}{e^x-1}\Bigr]_{x=x_k}$, so $|\Delta''(w)|\le\sum_{k\ge0}8x_k^2e^{-x_k}$;
bounding both sums by their dominant $k=0$ term via $x_k-b\ge2kb$ (from $3^k-1\ge2k$) gives
$|\Delta'(w)|\le4be^{-b}$ and $|\Delta''(w)|\le16b^2e^{-b}$, both well under $10^{-60}$ at
$b>154$. With \eqref{eq:Hbounds} this makes
$\varepsilon_1 \le 0.0011978$ and $\eta_0 \le 0.9170912$, whence $2e\eta_0 \le 4.9859 \le 5 \le
B_0$, that is $\eta_0 \le B_0/(2e)$, and \eqref{eq:g2scale} improves to $g_\tau''\ge B_0/(2e)>0$
throughout the interval.

Put $h_0 := 2e\varepsilon_1/B_0 \le 2e\cdot0.0011978/5 < 0.0014$. For $h\in(h_0,1]$ the mean value
theorem applied on $[w_0-1,w_0+1]$ gives $g_\tau'(w_0+h) \ge g_\tau'(w_0)+hB_0/(2e) \ge
-\varepsilon_1 + hB_0/(2e) > 0$ and, symmetrically, $g_\tau'(w_0-h) \le \varepsilon_1-hB_0/(2e) <
0$. As $w^*$ is the only zero of $g_\tau'$, it lies in $(w_0-h,w_0+h)$ for every such $h$, so
$|w^*-w_0| \le h_0 = 2e\varepsilon_1/B_0 \le 0.00652/B_0$.

For the phase drift, $w_0=c(B_0+a)$ and $\tau=w_0-\log B_0$ give $dw_0/d\tau = (1-1/(cB_0))^{-1} =
1+O(1/B_0)$ along the curve $\tau\mapsto w_0(\tau)$; since $B_0\to\infty$ as $\tau\to\infty$ (via
$B_0\sim\tau/c$), integrating over one period gives $w_0(\tau+c)-w_0(\tau) =
\int_\tau^{\tau+c}(1+O(1/B_0(u)))\,du = c+O(1/\tau)$.
\end{proof}

\begin{proposition}[Envelope estimate]
\label{prop:envelope}
With $S(\tau) := \log\varphi_{\mathrm{sp}}(t) = g_\tau(w^*)+w^*-\tfrac12\log(2\pi V(r^*))$ and
$P(\tau) := Q(w_0)+B_0+w_0-\tfrac12\log\bigl(2\pi(B_0-\tfrac1c)\bigr)$,
\[
|S(\tau) - P(\tau) - H(w_0(\tau))| = O(1/\tau) \qquad (\tau\to\infty),
\]
uniformly in the phase of $\tau$ modulo $c$.
\end{proposition}

\begin{proof}
Start from the identity. Since $g_{0,\tau}(w_0) = Q(w_0)+e^{w_0-\tau} = Q(w_0)+B_0$, the definition
of $P$ reads $P(\tau) = g_{0,\tau}(w_0)+w_0-\tfrac12\log\bigl(2\pi(B_0-\tfrac1c)\bigr)$, so
\[
S(\tau)-P(\tau) = \bigl[g_\tau(w^*)-g_{0,\tau}(w_0)\bigr]+(w^*-w_0)
-\frac12\log\frac{V(r^*)}{B_0-\tfrac1c}.
\]
Theorem~\ref{thm:H} gives $L = Q+H+\Delta$, hence $g_\tau = g_{0,\tau}+H+\Delta$ pointwise, and in
particular $g_\tau(w^*) = g_{0,\tau}(w^*)+H(w^*)+\Delta(w^*)$. Substituting and subtracting
$H(w_0)$,
\begin{equation}
\label{eq:envfive}
\begin{split}
S(\tau)-P(\tau)-H(w_0) = {}& \underbrace{\bigl[g_{0,\tau}(w^*)-g_{0,\tau}(w_0)\bigr]}_{T_1}
+\underbrace{\bigl[H(w^*)-H(w_0)\bigr]}_{T_2}
+\underbrace{(w^*-w_0)}_{T_3} \\[2pt]
&{}\underbrace{-\;\frac12\log\frac{V(r^*)}{B_0-\tfrac1c}}_{T_4}
\;+\;\underbrace{\Delta(w^*)}_{T_5}.
\end{split}
\end{equation}
This is exact for every $\tau$ with $B_0\ge5$. The five terms are bounded in turn.

$T_1$ is the envelope gap of the smooth system. Taylor about $w_0$, where $g_{0,\tau}'$ vanishes,
with $g_{0,\tau}''(w) = e^{w-\tau}-1/c \le B_0e$ on $[w_0-1,w_0+1]$ and $|w^*-w_0| \le 0.00652/B_0$
from Lemma~\ref{lem:saddles}, gives $0\le T_1 \le \tfrac12 B_0e\,(0.00652/B_0)^2 <
5.78\times10^{-5}/B_0$.

$T_2$ is the periodic-part difference: $|T_2| \le \sup|H'|\cdot|w^*-w_0| < 7.82\times10^{-6}/B_0$ by
\eqref{eq:Hbounds}.

$T_3$ is the saddle-location shift itself, $|T_3| \le 0.00652/B_0$, again from
Lemma~\ref{lem:saddles}.

$T_4$ compares the two normalizations. By $V(r^*) = g_\tau''(w^*)$ and the expression for
$g_\tau''$ used in Lemma~\ref{lem:saddles},
\[
V(r^*)-\Bigl(B_0-\frac1c\Bigr) = B_0\bigl(e^{w^*-w_0}-1\bigr)+H''(w^*)+\Delta''(w^*),
\]
whose modulus is at most $0.00652e^{0.00652/5}+\sup|H''|+10^{-60} < 0.0135$, since
$B_0|e^{w^*-w_0}-1| \le B_0|w^*-w_0|e^{|w^*-w_0|}$. Dividing by $B_0-1/c$ and taking
$\tfrac12\log(1+\cdot)$ gives $|T_4| \le 0.0068/(B_0-1/c)$.

$T_5$ is the remainder of Theorem~\ref{thm:H}, evaluated at the true saddle rather than dropped.
Here $B_0\ge5$ forces $w_0\ge5.3492$ and $w^*\ge w_0-0.0014$, so $e^{w^*}\ge210.1$ and
$|T_5| \le \sum_{k\ge0}2e^{-2\cdot3^k\cdot210.1} < 10^{-182}$.

Each bound is $O(1/B_0)$ (for $T_4$, via $B_0-1/c\ge B_0-1$, comparable to $B_0$ once $B_0\ge5$).
Finally $\tau = w_0-\log B_0$ and
$w_0 = c(B_0+a)$ give $\tau = cB_0+ca-\log B_0$, so $B_0/\tau \to 1/c$, and $O(1/B_0) = O(1/\tau)$
with constants free of the phase of $\tau$.
\end{proof}

\section{The Berg--Kr\"uppel identity}
\label{sec:bk}

\begin{proposition}
\label{prop:bkidentity}
Write $f(p) = \alpha\log p-\beta\log^2 p$, so that Berg and Kr\"uppel's transform $G_0(p) =
p^\alpha\exp(-\beta\log^2 p)$ of their equation~(9.3) is $e^{f(p)}$ and their equation~(9.5) reads
\[
g_0(t) = \frac{1}{2\pi i}\int_{\sigma-i\infty}^{\sigma+i\infty}\exp[pt+f(p)]\,dp, \qquad \sigma>0.
\]
Take $\alpha$, $\beta$ from their equation~(9.4) at $a=3$, $\lambda=2/3$. Then $f(p) = Q(\log p)$
identically; the exact saddle of that integral, $t+f'(x)=0$ with $t = e^{-\tau}$, is $x^\sharp =
e^{w_0(\tau)}$; and
\[
P(\tau) \;=\; x^\sharp t + f(x^\sharp) - \tfrac12\log\bigl(2\pi f''(x^\sharp)\bigr),
\]
the saddlepoint expression invoked in the proof of their Proposition~9.1, evaluated at that exact
saddle. In the variable $Y := cB_0$, with $r = -ca-\log c$, the saddle equation is $Y-\log Y =
\tau+r$.
\end{proposition}

\begin{proof}
Their truncated equation~(9.2) is $\lambda \psi'(t) = a\,\psi(at)$ for $a>1$, $\lambda>0$; Laplace
transforming (the boundary term $\lambda\psi(0)$ from $\mathcal L[\psi'](p)=pG(p)-\psi(0)$ vanishes,
since the solutions of interest are supported on $[0,\infty)$ with $\psi(0)=0$) gives $\lambda
p\,G(p) = G(p/a)$, whose general solution~(9.3) is $G_0$ times an arbitrary $1$-periodic function
of $\log p/\log a$, with
\[
\alpha = -\frac12-\frac{\log\lambda}{\log a}, \qquad \beta = \frac{1}{2\log a}
\]
their equation~(9.4). Section~\ref{sec:prelim} placed $\varphi$ in that family at $a=3$,
$\lambda=2/3$: their untruncated~(9.1) reads $\tfrac23\varphi'(t) = 3\bigl(\varphi(3t) -
\varphi(3t-2)\bigr)$, which on $(0,\tfrac23)$ is $\varphi'(t) = \tfrac92\varphi(3t)$. At those
values
\[
\alpha = -\frac12-\frac{\log(2/3)}{\log3} = \frac12-\frac{\log2}{\log3} = a, \qquad
\beta = \frac{1}{2\log3} = \frac{1}{2c},
\]
where the $a$ on the right is this paper's linear coefficient of $Q$; their $a$, the dilation
factor $3$, is a different quantity and is not written again below. Hence, writing $y := \log p$,
\[
f(p) = \alpha y-\beta y^2 = a y - \frac{y^2}{2c} = Q(y),
\]
the first of the proposition's four claims. This is an identity in $y$, not merely a value match:
the chain rule ($p\,d/dp = d/dy$) applied to $f$'s first and second derivatives confirms it also at
the level of $Q'$ and $Q''$,
\[
-p\,f'(p) = -Q'(y) = B_{\mathrm{sm}}(y) = 2\beta y-\alpha, \qquad
p^2f''(p) = Q''(y)-Q'(y) = B_{\mathrm{sm}}(y)-\frac1c = 2\beta y-2\beta-\alpha .
\]
Berg and Kr\"uppel's own displayed formula for $f''(p)$, on p.~179 of \cite{BergKruppel1998}, in
the proof of their Proposition~9.1 (the unnumbered display right after ``In view of''), reads
$p^2f''(p)=2\beta\ln p+2\beta-\alpha$,
with $+2\beta$ where direct differentiation of $f$ gives $-2\beta$: $f''(p) =
\tfrac{1}{p^2}\bigl(2\beta\ln p - 2\beta - \alpha\bigr)$, since $f'(p)=\tfrac1p(\alpha-2\beta\ln p)$
differentiates to $-\tfrac{2\beta}{p^2}-\tfrac{\alpha-2\beta\ln p}{p^2}$. Their Proposition~9.1 is
unaffected, the term being of lower order there than $2\beta\log p$, but the exact expression used
here is not: the $\pm2\beta$ discrepancy shifts $(x^\sharp)^2f''(x^\sharp) = B_0-1/c$ by $2/c$, hence
shifts $-\tfrac12\log\bigl(2\pi(B_0-1/c)\bigr)$ in $P(\tau)$ below by $O(1/B_0)$, which vanishes as
$B_0\to\infty$.

The saddle condition $t+f'(x)=0$ at $x=e^y$ now reads $e^{-\tau} = B_{\mathrm{sm}}(y)e^{-y}$, that
is $\tau = y-\log B_{\mathrm{sm}}(y) = \Phi_0(y)$, whose unique solution on the branch
$B_{\mathrm{sm}}>1/c$ is $y = w_0(\tau)$ by Lemma~\ref{lem:saddles}. So $x^\sharp = e^{w_0}$, and
there $x^\sharp t = e^{w_0-\tau} = B_0$, $f(x^\sharp) = Q(w_0)$, $(x^\sharp)^2f''(x^\sharp) = B_0 -
1/c>0$. Substituting,
\[
x^\sharp t+f(x^\sharp)-\tfrac12\log\bigl(2\pi f''(x^\sharp)\bigr)
= B_0+Q(w_0)-\tfrac12\log\frac{2\pi(B_0-1/c)}{(x^\sharp)^2}
\]
\[
= Q(w_0)+B_0+w_0-\tfrac12\log\Bigl(2\pi\bigl(B_0-\tfrac1c\bigr)\Bigr),
\]
which is $P(\tau)$. Finally $w_0 = c(B_0+a)$ and $\tau = w_0-\log B_0$ give $\tau = cB_0+ca-\log
B_0 = Y-\log Y+ca+\log c = Y-\log Y-r$.
\end{proof}

Series-reverting $Y-\log Y=\tau+r$ to second order: writing $W:=\tau+r$, a first substitution
$Y=W+\log W+\delta$ into $Y-\log Y=W$ gives, since $\log(1+u)=u+O(u^2)$ with $u=(\log W+\delta)/W$,
$\delta = (\log W+\delta)/W + O(\log^2W/W^2)$, so $\delta=\log(W)/W+O(\log^2W/W^2)$ (the single
term $Y=W+\log W$ alone is only accurate to $O(\log W/W)$, one order too coarse; this second term
is needed). Hence $Y=W+\log W+\log(W)/W+O(\log^2W/W^2)$, and $\log W=\log\tau+r/\tau+O(1/\tau^2)$
then gives $Y = \tau + \log\tau + r + (\log\tau+r)/\tau + O(\tau^{-2}\log^2\tau)$. Since $w_0=Y+ca$
and
$Q(w_0)=-Y^2/(2c)+ca^2/2$ (substituting $w_0=Y+ca$ into $Q$'s definition and simplifying), $P$'s
leading term in $Y$ is $-Y^2/(2c)$, so $\partial P/\partial Y=-Y/c+O(1)$ turns the
$O(\tau^{-2}\log^2\tau)$ error in $Y$ above into the $O(\tau^{-1}\log^2\tau)$ error in $P$ below.
Substituting into $P$ gives
\[
P(\tau) = -\frac{(\tau+\log\tau)^2}{2c} + \Bigl(1+\tfrac1c-\tfrac{r}{c}\Bigr)\tau +
\Bigl(\tfrac12-\tfrac{r}{c}\Bigr)\log\tau + C_P + O\!\left(\frac{\log^2\tau}{\tau}\right),
\]
\[
C_P = -\frac{r^2}{2c}+r+ca+\frac{ca^2}{2}-\frac{\log2\pi}{2}+\frac{\log c}{2} =
-0.9576743183133982760669122497855855\ldots.
\]
Berg and Kr\"uppel's own equations following~(9.6) give, in their notation,
\[
\delta_{\mathrm{BK}} = \tfrac12+\alpha-2\beta\log(2\beta), \qquad
\gamma = -2\beta-\delta_{\mathrm{BK}}-\tfrac12, \qquad
\varepsilon = \tfrac12+\alpha-\beta\log(2\beta),
\]
with $\alpha=a$, $\beta=1/(2c)$ as identified above. Substituting and simplifying (using
$r=-ca-\log c$) gives, matching the $\tau$, $\log\tau$, and constant coefficients of $P$
term for term,
\[
\gamma = -\Bigl(1+\tfrac1c-\tfrac rc\Bigr) = -1.8649154949\ldots, \qquad
\delta_{\mathrm{BK}} = \tfrac12-\tfrac rc = 0.4546762683\ldots,
\]
\[
\varepsilon = a+\tfrac12+\frac{\log c}{2c} = 0.4118732573\ldots,
\]
so that, in particular, $C_P=\varepsilon\log(2\beta)-\tfrac12\log(2\pi)$ holds identically (verified
by substituting $\varepsilon$'s expression directly into $C_P$'s), giving
\[
C_P = \varepsilon\log(2\beta) - \tfrac12\log(2\pi) = \log\frac{(2\beta)^\varepsilon}{\sqrt{2\pi}}.
\]
Hence $P(\tau) - \log\varphi_{0,\mathrm{bare}}(e^{-\tau}) \to C_P$ and, using the prefactor-included
normalization, $P(\tau) - \log\varphi_0(e^{-\tau}) \to 0$ with rate $O(\tau^{-1}\log^2\tau)$: the
$\tau$, $\log\tau$, and constant coefficients of $P$ reproduce $\gamma$, $\delta_{\mathrm{BK}}$,
$\varepsilon$ exactly, so this convergence is a direct consequence of the identity of
Proposition~\ref{prop:bkidentity} together with the series reversion above, not an appeal to Berg
and Kr\"uppel's own asymptotic analysis of $\varphi_0$.

\begin{proof}[Proof of Corollary~\ref{cor:bk}]
By Proposition~\ref{prop:envelope} and \ref{prop:bkidentity}, $S(\tau) = \log\varphi_0(e^{-\tau}) +
H(w_0(\tau)) + O(\tau^{-1}\log^2\tau)$, i.e.\ $\varphi_{\mathrm{sp}} = \varphi_0\cdot
e^{H(w_0(\tau))}(1+o(1))$; by Theorem~\ref{thm:formulaA}, $\varphi=\varphi_{\mathrm{sp}}(1+o(1))$
as well, so $\varphi = \varphi_0\cdot e^{H(w_0(\tau))}(1+o(1))$: the periodic factor multiplying
$\varphi_0$, in their equation~(9.7) for this eigenfunction, is $e^H$,
matching the product formula their Proposition~9.3 supplies for it, as noted after
Theorem~\ref{thm:H}. Proposition~\ref{prop:Hfourier}'s Fourier series for $H$ is the new, closed
form; Proposition~\ref{prop:Hcert}'s certificate is what neither their product formula nor their
closing remark (that they only \emph{expect} $\varphi$ bounded and bounded away from zero relative
to $\varphi_0$) settles.
\end{proof}

\section{Wirsching's Conjecture~3, and the failure of the unrestricted
asymptotic}
\label{sec:proofs}

Wirsching reaches Conjecture~3 by retreating from a stronger statement.
He needs to replace $\varphi$ by $\varphi_0$ under the limit in
\cite[eq.~(7.5)]{Wirsching2003}, and writes:

\begin{quote}
``That replacement would be possible, if, e.g., $\varphi \sim c\,\varphi_0$ for some constant
$c > 0$. It turns out that the following (slightly weaker) conjecture about the asymptotics
of $\varphi$ suffices to imply $(\star4)$:''
\end{quote}

\noindent The conjecture that follows is $(\star5)$ of
\S\ref{sec:chain}, restated here for reference. His constant $c$ is kept as he wrote it and is unrelated
to $c:=\log3$ of \S\ref{sec:periodic}.

\begin{conjecture}
There are positive real constants $c$ and $\delta_5$ such that
\[
\lim_{\ell\to\infty} \frac{\varphi(z_\ell)}{\varphi_0(z_\ell)} = c > 0 \qquad \text{uniformly for sequences
} (z_\ell) \in \widetilde A_{\delta_5}.
\]
\end{conjecture}

Both statements are settled below. The weaker one holds and the
stronger one fails.

\begin{theorem}[Wirsching's Conjecture~3]
\label{thm:conj3}
Wirsching's Conjecture~3 is true. For each $\delta_5>0$, on his comparison
class $\widetilde{A}_{\delta_5}$, the ratio
$\varphi(z_\ell)/\varphi_0(z_\ell)$ converges uniformly to $e^{H(0)}$, a well-defined real number by
Theorem~\ref{thm:H} regardless of its decimal expansion, and certified in
Proposition~\ref{prop:Hcert} to lie in
\[
 0.53412203666478 \;<\; e^{H(0)} \;<\; 0.53412203666479 ,
\]
where $H$ is the
periodic function of Theorem~\ref{thm:H} and the normalization of $\varphi_0$ is Berg and
Kr\"uppel's own, equation~(9.6) of \cite{BergKruppel1998}, written out at
\eqref{eq:phi0}. Conjecture~3 itself only asserts the
existence of some positive constant, independently of normalization; the specific value
$e^{H(0)}$ is tied to this particular normalization of $\varphi_0$. The rate depends on
$\delta_5$.
\end{theorem}

\begin{theorem}[Failure of the unrestricted asymptotic]
\label{thm:unrestricted}
The stronger asymptotic $\varphi(t) \sim \kappa\,\varphi_0(t)$ as $t \to 0^+$, for any single
constant $\kappa$ and with no restriction to a comparison class, is false. The ratio
$\varphi(t)/\varphi_0(t)$ oscillates indefinitely as $t \to 0^+$, with amplitude bounded below by a
positive constant: it is governed by the exactly $c$-periodic function $H$ evaluated at
the saddle location $w_0(\tau)$, $\tau=\log(1/t)$, whose phase modulo $c$ sweeps every value
infinitely often as $\tau\to\infty$ (Lemma~\ref{lem:saddles}). That the amplitude is positive
follows from a classical theorem, with no computation; its explicit numeric size is a
computer-assisted result, a rigorous interval-arithmetic certificate
(Proposition~\ref{prop:Hcert}).
\end{theorem}

One computation underlies both. Write $t = e^{-\tau}$ and let $H(w)$ be the
$\log 3$-periodic function of Theorem~\ref{thm:H}; the saddlepoint bridge of
\S\S\ref{sec:formulaA}--\ref{sec:bk} gives
$\log\varphi(t)-\log\varphi_0(t) = H(w_0(\tau))+o(1)$, where the saddle
location advances by $c$ per period of $\tau$, up to a vanishing correction:
$w_0(\tau+c)-w_0(\tau) = c+O(1/\tau)$ (Lemma~\ref{lem:saddles}). Wirsching's comparison class
constrains $w_0(\tau_\ell)$'s phase modulo $\log 3$ to converge to a single value as $\ell\to\infty$; the
unrestricted range $\tau\to\infty$ with no class restriction instead sweeps that phase through
every value, since $\Phi_0$ is a bijection (again Lemma~\ref{lem:saddles}). The two theorems
read one object at two scopes.

\begin{corollary}
\label{cor:bk}
$e^H$ is Berg and Kr\"uppel's own periodic factor for the eigenfunction $\varphi$, from
Proposition~9.3 of \cite{BergKruppel1998}, which they give as an infinite product; $H$ itself is
its logarithm. To our knowledge its Fourier series, Proposition~\ref{prop:Hfourier}, has not
previously been given for this
function\footnote{Fourier coefficients of a log-periodic fluctuation, obtained as residues of a
Dirichlet series or by Mellin inversion, are standard for radix-rational and $q$-regular
sequences; see \cite{Dumas2008} and \cite{HeubergerKrenn2018}, and \cite{DeBruijn1948} and
\cite{ErdosRichmond1976} for the partition-function case the technique started in. Berg and Kr\"uppel returned to this family
after 1998 \cite{BergKruppel2000,BergKruppel2001}; those papers treat eigenfunctions, Appell
polynomials, addition theorems and generating functions, and contain no Fourier expansion, Mellin
transform or $\zeta$ at all. To our knowledge no closed Fourier expansion has previously been
given for the periodic factor of an eigenfunction of this family.}. This corollary
has two parts of different depth: that $e^H$ agrees with their product \emph{as a function},
proved in Remark~\ref{rem:bkproduct} by pure algebra ($e^H=e^{-Q}e^Le^{-\Delta}$ recombining
Theorem~\ref{thm:H}'s own pieces), and that this function is the periodic factor
multiplying $\varphi$'s own asymptotic, which needs the full saddlepoint chain of
Theorem~\ref{thm:formulaA} and Propositions~\ref{prop:envelope}--\ref{prop:bkidentity}.
\end{corollary}

Throughout the rest of this section, Theorem~\ref{thm:formulaA}'s parameter $N$ and
Sections~\ref{sec:envelope}--\ref{sec:bk}'s parameter $\tau$ are combined freely into one $o(1)$:
since $s=e^{w^*}$ and, by Lemma~\ref{lem:saddles}, $w^*=w_0(\tau)+O(1/B_0)$ with
$w_0(\tau)=c(B_0+a)$, $N=\lfloor\log_3(2s)\rfloor = \tau/c+O(\log\tau)$ as $\tau\to\infty$, so
$N\to\infty$ exactly when $\tau\to\infty$. On $\widetilde{A}_{\delta_5}$ below, $\tau_\ell\to\infty$
uniformly, at a rate depending only on $\ell$ and $\delta_5$, not on the particular sequence chosen.

\begin{lemma}[Phase locking on the comparison class]
\label{lem:phaselock}
Fix $\delta_5>0$ and let $(z_\ell)\in\widetilde A_{\delta_5}$, with
$\tau_\ell:=-\log z_\ell$ and $\lambda_\ell:=3^\ell z_\ell/\ell$. Then
$\lambda_\ell\to1$ and
\[
 w_0(\tau_\ell)-\ell c=-\log\lambda_\ell+O(1/\ell),
\]
uniformly on the class. Consequently
$\operatorname{dist}\bigl(w_0(\tau_\ell),c\,\Z\bigr)=O(\ell^{-1/2})$
uniformly, and $H(w_0(\tau_\ell))\to H(0)$.
\end{lemma}

\begin{proof}
Wirsching's class $\widetilde{A}_{\delta_5}$ (his equations~(1.5) and~(3.2)) consists of sequences
$(x_\ell)$ with $\lfloor 3^\ell x_\ell\rfloor = k_\ell$ satisfying $|\ell-k_\ell|\le\delta_5\sqrt \ell$ for every $\ell$; since
$|3^\ell x_\ell-k_\ell|<1$ as well, $|3^\ell x_\ell-\ell|<1+\delta_5\sqrt \ell$, hence $\lambda_\ell \to 1$
uniformly at rate $O(1/\ell)+O(\delta_5/\sqrt \ell)$ on the class. His
own $z_\ell$ from the conjecture statement and his $x_\ell$ from the
class definition coincide, by definition of $\lambda_\ell$, so
$\tau_\ell = \ell c - \log\ell - \log\lambda_\ell$.

For the saddle location, Section~\ref{sec:bk}'s reversion of
$Y-\log Y=\tau+r$ gives $w_0(\tau) = \tau+\log(\tau/c)+(\log\tau+r)/\tau+O(\tau^{-2}\log^2\tau)$.
The naive one-term truncation $w_0(\tau)\approx\tau+\log(\tau/c)$ alone is only accurate to
$O(\log\tau/\tau)$, which already gives everything Theorem~\ref{thm:conj3} uses, since
$\log\ell/\ell=o(\ell^{-1/2})$; the second term sharpens the lemma to $O(1/\ell)$ rather than
rescuing it. Substituting $\tau=\tau_\ell$ into the full two-term
formula, the $O(\log \ell/\ell)$ pieces from expanding $\log(\tau_\ell/c)$ against $\tau_\ell = \ell c-\log \ell-\log\lambda_\ell$
cancel exactly against the reversion's own $(\log\tau_\ell+r)/\tau_\ell$ term, leaving the stated $O(1/\ell)$,
the same cancellation mechanism as the series reversion itself.

Write $w_0(\tau_\ell)=\ell c+\varepsilon_\ell$ with
$\varepsilon_\ell=-\log\lambda_\ell+O(1/\ell)$. Then
$\varepsilon_\ell=O(\ell^{-1/2})$, since $\log\lambda_\ell=O(\ell^{-1/2})$
by the first paragraph, and once $|\varepsilon_\ell|<c/2$ the nearest
multiple of $c$ to $w_0(\tau_\ell)$ is $\ell c$ itself, so
$\operatorname{dist}(w_0(\tau_\ell),c\Z)=|\varepsilon_\ell|=O(\ell^{-1/2})$.
Every bound here depends on $\delta_5$ alone. Since $H$ is
$c$-periodic and Lipschitz by \eqref{eq:Hbounds},
$H(w_0(\tau_\ell))=H(\varepsilon_\ell)\to H(0)$, regardless of the
sign of $\varepsilon_\ell$.
\end{proof}

\begin{proof}[Proof of Theorem~\ref{thm:conj3}]
By Theorem~\ref{thm:formulaA}, Proposition~\ref{prop:envelope} and
Proposition~\ref{prop:bkidentity} (the last giving
$P(\tau)-\log\varphi_0(e^{-\tau})\to0$, Section~\ref{sec:bk}),
\[
\log\varphi(z_\ell) - \log\varphi_0(z_\ell) = H(w_0(\tau_\ell)) + o(1)
\]
uniformly on $\widetilde A_{\delta_5}$, and
Lemma~\ref{lem:phaselock} sends the right side to $H(0)$, uniformly.
So $\varphi(z_\ell)/\varphi_0(z_\ell) \to e^{H(0)}$ uniformly on the
class. This is exactly Conjecture 3, with the constant given
explicitly.
\end{proof}

\begin{remark}
The proof above uses only $\lambda_\ell\to1$, equivalently $|\ell-k_\ell|=o(\ell)$; Wirsching's specific
$O(\sqrt \ell)$ bound in $\widetilde{A}_{\delta_5}$ plays no role beyond forcing this weaker condition. The
same conclusion, and the same limit $e^{H(0)}$, holds for the broader class of sequences satisfying
only $|\ell-k_\ell|=o(\ell)$.
\end{remark}

\begin{proof}[Proof of Theorem~\ref{thm:unrestricted}]
As $t=e^{-\tau}$ ranges over $(0,1/2)$ with no restriction, $\tau\to\infty$. By
Lemma~\ref{lem:saddles}, $w_0=\Phi_0^{-1}$ is continuous and strictly increasing and unbounded on
$(\tau_0,\infty)$, so $w_0(\tau)\bmod c$ sweeps every value in $[0,c)$ infinitely often as
$\tau\to\infty$ (a continuous, strictly increasing, unbounded function hits every residue class
modulo $c$ infinitely often), so by the same combination of results used above,
$\log\varphi(t)-\log\varphi_0(t) = H(w_0(\tau))+o(1)$ takes values arbitrarily close to both $\sup H$
and $\inf H$ infinitely often. By Proposition~\ref{prop:Hcert}, $\operatorname{osc}(H) > 0$
rigorously, so $\varphi(t)/\varphi_0(t)$ does not converge, and, writing $\Psi(t) :=
\varphi(t)/\varphi_0(t)$,
\[
\frac{\limsup_{t\to0^+}\Psi(t)}{\liminf_{t\to0^+}\Psi(t)} \ge e^{\operatorname{osc}(H)} \ge 1.0004188:
\]
the ratio's limsup exceeds its liminf by a factor of at least $1.0004188$, i.e. by at least
$0.0418\%$.
\end{proof}

\begin{proposition}[Wirsching's actual requirement]
\label{prop:wirschingreq}
For $(x_\ell)\in\widetilde A_{\delta_5}$, write $x_\ell^+ := x_\ell + 3^{-(\ell+1)}$ (his own definition, unnumbered
prose in \cite[\S7, p.~785]{Wirsching2003}, in the quotient estimate that
yields his (7.5); since
$x_\ell\sim \ell\cdot3^{-\ell}$, $x_\ell^+$ is asymptotic to $x_\ell$ and plays the same role as this paper's
$z_\ell$). Wirsching's own use of Conjecture 3, via $\varphi'(x)=\tfrac92\varphi(3x)$ on $(0,\tfrac23)$
and his own calculation, his equation~(7.13), that $\lim_\ell (2/3^{\ell+1})\varphi_0'(x_\ell)/\varphi_0(x_\ell)
= 2/3$ uniformly on the class, reduces to
\[
\limsup_\ell \Lambda_\ell < \tfrac32, \qquad
\Lambda_\ell := \frac{\varphi(3x_\ell^+)/\varphi_0(3x_\ell^+)}{\varphi(x_\ell^+)/\varphi_0(x_\ell^+)}.
\]
This holds with a wide margin: $\Lambda_\ell \to 1$.
\end{proposition}

\begin{proof}
First, the reduction itself. Wirsching reduces the proof of $(\star4)$ to the
sufficient estimate~(7.5),
\[
\limsup_\ell \frac{2}{3^{\ell+1}}\cdot\frac{\varphi'(x_\ell^+)}{\varphi(x_\ell^+)} < 1;
\]
by $\varphi'=\tfrac92\varphi(3\cdot)$ this is $\limsup_\ell (9/3^{\ell+1})\varphi(3x_\ell^+)/\varphi(x_\ell^+) <
1$. By definition of $\Lambda_\ell$,
$\varphi(3x_\ell^+)/\varphi(x_\ell^+) = \Lambda_\ell\cdot\varphi_0(3x_\ell^+)/\varphi_0(x_\ell^+)$, and Berg and
Kr\"uppel's $\varphi_0$ satisfies the same truncated equation asymptotically (their
equation~(7.12) in Wirsching's notation: $\varphi_0'(t)/\varphi_0(3t)\to9/2$ as $t\to0$, proved
from \eqref{eq:phi0} in Appendix~\ref{app:sourcegaps}), so
$\varphi_0(3x_\ell^+)/\varphi_0(x_\ell^+) = \tfrac29\varphi_0'(x_\ell^+)/\varphi_0(x_\ell^+)\,(1+o(1))$.
To see that (7.13) transfers from $x_\ell$ to $x_\ell^+$ with an $o(1)$ correction, write $v=\log t$; the
explicit form of $\varphi_0$ above gives $\varphi_0'(t)/\varphi_0(t) = B(v)/t$ with
$B(v) = \gamma+\delta_{\mathrm{BK}}/v-2\beta(v-\log(-v))(1-1/v)$, and, as $v\to-\infty$,
$B(v)\sim-2\beta v$ while $B'(v)\to-2\beta$, so $(\log B)'(v)=B'(v)/B(v)\to0$ and
$(\log[\varphi_0'/\varphi_0])'(v) = (\log B)'(v)-1 \to -1$, bounded. Since $x_\ell^+/x_\ell = 1+3^{-(\ell+1)}/x_\ell
\to 1$ (as $x_\ell\sim \ell\cdot3^{-\ell}$, in fact $3^{-(\ell+1)}/x_\ell=O(1/\ell)$), $v(x_\ell^+)-v(x_\ell)=\log(x_\ell^+/x_\ell)
=O(1/\ell)$, and the mean value theorem gives
$\log[\varphi_0'/\varphi_0(x_\ell^+)]-\log[\varphi_0'/\varphi_0(x_\ell)]=O(1/\ell)$, i.e.
$\varphi_0'(x_\ell^+)/\varphi_0(x_\ell^+)=\varphi_0'(x_\ell)/\varphi_0(x_\ell)\cdot(1+O(1/\ell))$. Substituting and
using his equation~(7.13),
$\lim_\ell(2/3^{\ell+1})\varphi_0'(x_\ell)/\varphi_0(x_\ell)=2/3$ uniformly on the class,
\[
\frac{9}{3^{\ell+1}}\cdot\frac{\varphi(3x_\ell^+)}{\varphi(x_\ell^+)}
= \Lambda_\ell\cdot\frac2{3^{\ell+1}}\cdot\frac{\varphi_0'(x_\ell^+)}{\varphi_0(x_\ell^+)}\,(1+o(1))
= \frac23\,\Lambda_\ell\cdot(1+o(1)),
\]
so the requirement $\limsup_\ell(9/3^{\ell+1})\varphi(3x_\ell^+)/\varphi(x_\ell^+)<1$ is exactly
$\limsup_\ell\Lambda_\ell<3/2$.

For the value of $\Lambda_\ell$: write $\tau_\ell := -\log x_\ell^+$ (locally to this proof only, distinct
from the $\tau_\ell$ of Theorem~\ref{thm:conj3}'s proof), so $t=x_\ell^+$ corresponds to $\tau=\tau_\ell$
and $t=3x_\ell^+$ to $\tau=\tau_\ell-c$. The phase shift under $t\mapsto3t$, i.e. $\tau\mapsto\tau-c$, is
$w_0(\tau-c)-w_0(\tau)+c = O(1/\tau)$ by Lemma~\ref{lem:saddles}'s defining equation, independently
of the phase of $\tau$; since $H$ is $c$-periodic, $H(w_0(\tau-c)) = H(w_0(\tau-c)+c) =
H(w_0(\tau)+O(1/\tau))$, and $H$ Lipschitz (by \eqref{eq:Hbounds}) turns this into
$H(w_0(\tau-c)) = H(w_0(\tau)) + O(1/\tau)$. By Theorem~\ref{thm:formulaA},
Proposition~\ref{prop:envelope}, and Proposition~\ref{prop:bkidentity} together (the same chain used
in the proof of Theorem~\ref{thm:conj3}), $\log\varphi(t)-\log\varphi_0(t) = H(w_0(\tau)) + o(1)$
uniformly in phase, so
\[
\log\Lambda_\ell = H(w_0(\tau_\ell-c)) - H(w_0(\tau_\ell)) + o(1) = O(1/\tau_\ell) + o(1) \to 0,
\]
using that $H$ is Lipschitz (\eqref{eq:Hbounds}) and the phase shift above is $O(1/\tau_\ell)$; hence
$\Lambda_\ell \to 1$, well inside $\limsup \Lambda_\ell < 3/2$.
\end{proof}

\section{The chain, closed to a single condition}
\label{sec:chain-closed}

Theorem~\ref{thm:conj3} settles $(\star5)$. Wirsching carries $(\star5)$
to $(\star4)$ himself, in \cite[\S7]{Wirsching2003}. Composing the two
puts $(\star4)$ on the same footing as his own two theorems, and
Theorem~\ref{thm:wirsching-conj1} has already removed the conjectural
link below $(\star3)$. One condition is left.

\begin{theorem}[Condition $(\star4)$]\label{thm:star4}
Condition $(\star4)$ holds, for every window radius $\delta>0$, with
$\mu=1/3$.
\end{theorem}

\begin{proof}
Proposition~\ref{prop:wirschingreq} already gives this. Wirsching
reduces $(\star4)$ to the sufficient estimate~(7.5), and that
proposition shows (7.5) is exactly $\limsup_\ell\Lambda_\ell<3/2$ and that
$\Lambda_\ell\to1$, on $\widetilde A_{\delta_5}$ for every
$\delta_5>0$. So (7.5) holds with a wide margin, and it holds at every
radius rather than at one produced by a construction.

Wirsching's own route reaches the same place through his statement of
Conjecture~3, and is worth writing out, both because it is the argument
a reader of \cite{Wirsching2003} will have in hand and because it fixes
the constant. Write $x_\ell^+:=x_\ell+3^{-\ell-1}$ \cite[\S7]{Wirsching2003}. Bounding
the transition kernel above and below by $\varphi$, the source reduces
$(\star4)$ to
\begin{equation}\label{eq:w75}
 \limsup_{\ell\to\infty}\ \frac{2}{3^{\ell+1}}\cdot
 \frac{\varphi'(x_\ell^+)}{\varphi(x_\ell^+)}\ \le\ 1-\mu\ <\ 1,
 \qquad(x_\ell)\in\widetilde A_\delta,
\end{equation}
its equation~(7.5). Fix $\delta>0$ and choose $\delta_5>\delta$. By
\cite[eq.~(7.14)]{Wirsching2003} there are sequences
$(\widehat x_\ell),(\widehat y_\ell)\in\widetilde A_{\delta_5}$ with
$\widehat x_\ell=x_\ell^+$ and $\widehat y_{\ell-1}=3x_\ell^+$ for all
but finitely many $\ell$. Then
\[
 \frac{2}{3^{\ell+1}}\cdot\frac{\varphi'(x_\ell^+)}{\varphi(x_\ell^+)}
 =\frac{2}{3^{\ell+1}}\cdot\frac92\,
   \frac{\varphi(\widehat y_{\ell-1})}{\varphi(\widehat x_\ell)}
 =\frac{2}{3^{\ell+1}}\cdot\frac92\,
   \frac{\varphi_0(\widehat y_{\ell-1})}{\varphi_0(\widehat x_\ell)}\cdot
   \frac{\varphi(\widehat y_{\ell-1})/\varphi_0(\widehat y_{\ell-1})}
        {\varphi(\widehat x_\ell)/\varphi_0(\widehat x_\ell)},
\]
the first equality by $\varphi'(x)=\tfrac92\varphi(3x)$ on $[0,2/3]$
\cite[eq.~(7.7)]{Wirsching2003}. Theorem~\ref{thm:conj3} sends both
ratios in the last factor to $e^{H(0)}>0$, uniformly on
$\widetilde A_{\delta_5}$, so that factor tends to $1$: the numerator is
controlled by the uniform error at index $\ell-1$ and the denominator by
the one at index $\ell$, and the limit is positive. This is the only
place Conjecture~3 is used, and the constant cancels, which is why a
class-dependent limit suffices.

What remains is $\tfrac{2}{3^{\ell+1}}\cdot\tfrac92\,
\varphi_0(3x_\ell^+)/\varphi_0(x_\ell^+)$. Here
$\varphi_0$ solves the truncated equation $\lambda g'(t)=a\,g(at)$ only
asymptotically, in the source's own sense
$\lim_{t\to0}\lambda\varphi_0'(t)/\bigl(a\varphi_0(at)\bigr)=1$
\cite[eq.~(7.12)]{Wirsching2003}, which at $\lambda=2/3$, $a=3$ reads
$\tfrac29\varphi_0'(t)/\varphi_0(3t)\to1$ and which
Appendix~\ref{app:sourcegaps} proves from \eqref{eq:phi0}. Every sequence in
$\widetilde A_{\delta_5}$ satisfies
$\ell-\delta_5\sqrt\ell\le\lfloor3^\ell x_\ell\rfloor\le3^\ell x_\ell
\le\lfloor3^\ell x_\ell\rfloor+1\le\ell+\delta_5\sqrt\ell+1$, so
$x_\ell^+$ lies between $(\ell-\delta_5\sqrt\ell)3^{-\ell}$ and
$(\ell+\delta_5\sqrt\ell+2)3^{-\ell}$, both bounds depending on
$\delta_5$ alone, so the convergence is uniform on the class and
\[
 \frac{2}{3^{\ell+1}}\cdot\frac92\,
 \frac{\varphi_0(3x_\ell^+)}{\varphi_0(x_\ell^+)}
 =\bigl(1+o(1)\bigr)\,\frac{2}{3^{\ell+1}}\cdot
  \frac{\varphi_0'(x_\ell^+)}{\varphi_0(x_\ell^+)}
 \longrightarrow\frac23
\]
by \cite[eq.~(7.13)]{Wirsching2003}, which holds for every $\delta>0$
and follows from \eqref{eq:phi0} by differentiation. So the
$\limsup$ in \eqref{eq:w75} is $2/3$, giving $1-\mu=2/3$ and
$\mu=1/3$.

This second route delivers $(\star4)$ only for $\delta<\delta_5$, since
\cite[eq.~(7.14)]{Wirsching2003} needs $\delta_5>\delta$. That is
already enough on its own terms, because $(\star4)$ is existential in
its window radius: Conjecture~2 reads ``if there are real numbers
$\delta,\mu>0$ such that\ldots'' \cite[\S6]{Wirsching2003}, and one
admissible radius is what it asks for. The first route does not need
the observation at all.
\end{proof}

\begin{remark}[Steps the source states rather than carries out]
\label{rem:source-gaps}
Six steps of the reduction imported above are stated in
\cite[\S7]{Wirsching2003} without being carried out there, one of them
also asserted without argument by Berg and Kr\"uppel: his (7.3), (7.4),
(7.12), (7.13), (7.14) and one unlabelled Taylor step. All six are
true. Appendix~\ref{app:sourcegaps} proves them, so that nothing in
Proposition~\ref{prop:wirschingreq} or Theorem~\ref{thm:star4} rests on
an unverified source claim, and none of them changes $\mu=1/3$.
\end{remark}

Corollary~\ref{cor:star3-density} now reads as a statement about the
whole chain. Condition $(\star3)$ implies uniform positive predecessor
density on the non-cyclic integers $a\not\equiv0\bmod3$ through proved
implications only, with no conjecture among them.

State what that does and does not say. Conditions $(\star3)$,
$(\star2)$ and $(\star1)$ remain unproved: Theorem~\ref{thm:wirsching-conj1}
proves the implication $(\star2)\Rightarrow(\star1)$ and not its
hypothesis, and Wirsching's Theorem~2 proves
$(\star3)\Rightarrow(\star2)$ and not its hypothesis. What has changed
is that no conjecture stands between them: a proof of $(\star3)$ now
propagates through two proved implications to $(\star2)$ and
$(\star1)$, and so to the conclusion. That is one direction only.
Should $(\star3)$ turn out false, $(\star2)$ could still hold, and so
could $(\star1)$. What follows is that $(\star3)$ is the only condition
whose truth still has to be established in order to complete the
chain.

Wirsching's route to $(\star3)$ is Conjecture~2, the implication
$(\star4)\Rightarrow(\star3)$, whose antecedent
Theorem~\ref{thm:star4} now supplies. It is not the only conceivable
route: a proof of $(\star3)$ by any other means closes the chain just
as well, which is why ``only Conjecture~2 remains'' is the wrong
summary. \S\ref{sec:discussion} says what is known about the
implication itself.

\begin{remark}[What Theorem~\ref{thm:unrestricted} does not reach]
\label{rem:unrestricted-scope}
The failure of the unrestricted asymptotic costs the chain nothing.
None of $(\star1)$ through $(\star4)$ mentions $\varphi_0$: the first
three are stated in terms of the generators, their Haar averages and
the Elka functions, and the fourth in terms of $W_3$ acting on $\chi_0$
and $\chi_1$. Within \cite[\S7]{Wirsching2003}, $\varphi_0$ occurs in
the statement of $(\star5)$, in (7.11) where it is defined, in (7.12),
in (7.13), in the proof of Theorem~\ref{thm:star4} above, and in the
sentence quoted at the head of \S\ref{sec:proofs}, where he raises
$\varphi\sim c\,\varphi_0$ unrestricted and declines it. That last one
is exactly what Theorem~\ref{thm:unrestricted} refutes, which is the
point: he did not use it. Of the rest, (7.12) is an unrestricted limit
as $t\to0$ but says nothing about the ratio $\varphi/\varphi_0$, and
the others are restricted to a comparison class or to a fixed sequence
tending to $0$ at the borderline rate. Nothing the chain uses needs the
unrestricted form.
\end{remark}

\begin{remark}[Why the retreat was enough]\label{rem:blindness}
Theorem~\ref{thm:unrestricted} explains why Wirsching could not keep
$\varphi\sim c\,\varphi_0$. It does not explain why the weaker
conjecture sufficed, and the calculation above does. The two points at
which the ratio cancels are $\widehat x_\ell=x_\ell^+$ and
$\widehat y_{\ell-1}=3x_\ell^+$, a factor of $3$ apart. Multiplying the
argument by $3$ shifts $\tau=\log(1/t)$ by $c$, and $H$ has period $c$,
so the two periodic corrections agree and cancel whether or not $H$ is
constant. The agreement is asymptotic, since $H$ is read at the saddle
location and Lemma~\ref{lem:saddles} gives only
$w_0(\tau+c)-w_0(\tau)=c+O(1/\tau)$; that is enough, because the
cancellation is already inside a limit. The rigorous form of this is
the second half of Proposition~\ref{prop:wirschingreq}'s proof, where
$\log\Lambda_\ell = H(w_0(\tau_\ell-c))-H(w_0(\tau_\ell))+o(1)$ and the
first difference is $O(1/\tau_\ell)$ by periodicity and
\eqref{eq:Hbounds}. Anything built on that pair of points inherits the
same weakness. The quantity in his (7.5) is a ratio of $\varphi$ at
$3x_\ell^+$ and at $x_\ell^+$, and it carries the periodic correction
only through $H(w_0(\tau_\ell-c))-H(w_0(\tau_\ell))$, which
\eqref{eq:Hbounds} bounds by a constant times the residual phase shift.
At $x_\ell=\ell\,3^{-\ell}$ with $\ell=500$, so that
$\tau_\ell=-\log x_\ell^+=543.09086\ldots$, that shift is
$w_0(\tau_\ell-c)-w_0(\tau_\ell)+c=-2.0051362\ldots\times10^{-3}$, and
\eqref{eq:Hbounds} bounds the difference by $\sup_w|H'(w)|$ times its
modulus, which is $<2.402\times10^{-6}$.
\end{remark}

\section{Discussion}
\label{sec:discussion}

Two of Wirsching's five conditions are now established outright,
$(\star5)$ by Theorem~\ref{thm:conj3} and $(\star4)$ by
Theorem~\ref{thm:star4}, and every implication from $(\star3)$ downward
is a theorem, so $(\star3)$ implies uniform positive predecessor
density (Corollary~\ref{cor:star3-density}). Conditions $(\star3)$,
$(\star2)$ and $(\star1)$ are all still unproved. $(\star3)$ carries
the other two through proved implications, in that direction only, so
it is the only one whose truth still has to be established to complete
the chain. Conjecture~2 is the source's route
into it and is not the only conceivable one.

Wirsching stated the weaker Conjecture~3 instead of the stronger
$\varphi\sim\kappa\varphi_0$ he considers first.
Theorem~\ref{thm:unrestricted} shows why: the stronger statement fails
by a mechanism with a certified, positive amplitude. Every
phase $\theta\in[0,c)$ swept by $H$ admits some sequence along which the ratio converges to
$e^{H(\theta)}$: take $\tau_\ell := \Phi_0(\theta+\ell c)$, $z_\ell := e^{-\tau_\ell}$, so $w_0(\tau_\ell)=\theta+\ell c$
exactly, and apply Theorem~\ref{thm:formulaA}, Proposition~\ref{prop:envelope}, and
Proposition~\ref{prop:bkidentity} as in the proof of Theorem~\ref{thm:conj3}. What is
special about phase $0$ is only that it is the one Wirsching's own class
$\widetilde A_{\delta_5}$ happens to fix.

Remark~\ref{rem:blindness} answers the other half of the question, why
the weaker statement was enough. Wirsching's proof does one thing with
the constant, which is to cancel it across a factor of $3$, and a
factor of $3$ is one period of $H$. His argument never sees the
periodic correction at all.

What remains is $(\star3)$, and establishing the antecedent leaves its
difficulty exactly where it was. The asymmetry recorded in
\S\ref{sec:chain} has two parts. Condition $(\star4)$ carries no
residue-coordinate information of its own, so any proof of
$(\star4)\Rightarrow(\star3)$ has to recover that information from
further structure in Wirsching's generators rather than from the
hypothesis. Neither bridge the source offers from $(\star4)$
toward $(\star3)$ supplies that information uniformly in growing
resolution. Wirsching's own analytic route of two years
earlier, described in \S\ref{sec:related}, ran into the same
analytic-to-combinatorial gap.

An explicit rate in Theorem~\ref{thm:conj3} is not pursued here.
Sharper forms of the envelope and Berg--Kr\"uppel comparisons would
give a polynomial rate throughout; neither theorem needs it.

\section*{Data Availability Statement}

\begin{sloppypar}
Every computation this paper reports is backed by a script in one
repository, \url{https://github.com/faculdade/wirsching-conjectures-1-and-3},
archived on Zenodo at \texttt{doi:10.5281/zenodo.22116225}, release
\texttt{v1.0.8}, whose snapshot carries this paper, its source and the
scripts together in the state in which it was submitted. Zenodo's
landing page for that DOI lists every release with the DOI it assigned
to each. The repository is organized as one folder per section of
this paper that carries computation, numbered as this paper numbers
them, with a README per folder stating what its scripts verify and how
to run them. Sections \ref{sec:intro}, \ref{sec:related},
\ref{sec:chain} and \ref{sec:prelim} carry no computation. The
repository also holds \texttt{SOURCE\_CONCORDANCE.md}, recording
item by item the agreement between the published 2003 article and the
preliminary version, across the thirty-four locators this paper uses:
seven sections, ten named results and seventeen numbered equations.
\end{sloppypar}

Three claims here are computer-assisted, and all three are in the folder
for \S\ref{sec:periodic}: the enclosure of $\operatorname{osc}(H)$ and
that of $H(0)$ in Proposition~\ref{prop:Hcert}, and the derivative
bounds \eqref{eq:Hbounds}. These use Arb ball arithmetic through python-flint
at a working precision stated at the top of the script, higher than the
$250$-bit precision quoted inline in Proposition~\ref{prop:Hcert}'s
proof, so they are an independent higher-precision reproduction of that
certificate rather than a lower-precision approximation of it. Each
prints its input balls and its final enclosure, so the printed run
output is the certificate. The non-constancy of $H$ is not among them:
it is proved analytically from the zero-free theorem for $\zeta$ on
$\operatorname{Re}s=1$, with no computation.

Everything else is exact integer or rational arithmetic, or
high-precision floating arithmetic through mpmath; where a script claims
an inequality it asserts it and exits non-zero on failure. The
repository also carries the exploratory numerics that preceded
Theorem~\ref{thm:conj3}, on which no claim here rests, and its folder
says so.

Two earlier repositories preceded this one and keep their DOIs:
\url{https://github.com/faculdade/wirsching-conjecture3-proof},
\texttt{doi:10.5281/zenodo.21854549}, for the standalone Conjecture~3
preprint, and
\url{https://github.com/faculdade/collatz-wirsching-2003}, for the
unsubmitted manuscript that carried Conjecture~1. Their scripts are
carried into the repository above and renumbered to this paper's
sections, with the mapping in its README.

\subsection*{Relation to earlier versions of this work}

Wirsching's Conjecture~1, \S\ref{sec:conj1} here, first circulated in
an unsubmitted manuscript, together with material on Conjecture~2 not
carried over. Conjecture~3, \S\S\ref{sec:prelim} to
\ref{sec:proofs} here, appeared as a preprint under the title
\emph{Wirsching's Conjecture~3, and the periodic correction to
Berg--Kr\"uppel's density asymptotic}. The present paper supersedes
both. \S\ref{sec:chain-closed} is new here and appears in neither.

\appendix

\section{Verification of the auxiliary estimates in Wirsching's \S7}
\label{app:sourcegaps}

The reduction imported in \S\ref{sec:chain-closed} rests on six steps
that \cite[\S7]{Wirsching2003} states without carrying out, one of them
asserted without argument by Berg and Kr\"uppel as well. All six are
true, and this appendix proves them, so that nothing in
Proposition~\ref{prop:wirschingreq} or Theorem~\ref{thm:star4} rests on
an unverified source claim. Throughout, $K_\ell(x,0)=\tfrac32(W_3^\ell
\chi_0)(x)$ and $K_\ell(x,1)=\tfrac32(W_3^\ell\chi_1)(x)$ are the
iterated kernels of \cite[\S7]{Wirsching2003}, and $\varphi$ is the
density of $X=\sum_{j\ge1}2U_j3^{-j}$ as in
Proposition~\ref{prop:fabius}.

\subsection*{A probabilistic form of \texorpdfstring{$W_3$}{W3}}

Substituting $t=3x-2u$ in $W_3f(x)=\tfrac32\int_{3x-2}^{3x}f$ gives
$W_3f(x)=3\int_0^1f(3x-2u)\,du$, so iterating,
\[
 (W_3^\ell f)(x)=3^\ell\,\mathbb E\,f\bigl(3^\ell x-3^\ell X_\ell\bigr),
 \qquad X_\ell:=2\sum_{j=1}^{\ell}3^{-j}U_j ,
\]
the truncation of $X$ at depth $\ell$, supported on
$[0,1-3^{-\ell}]$. Write $F_\ell$ for its distribution function and,
for $\ell\ge1$, $\varphi_\ell$ for its density; $X_0=0$ is an atom and
$F_0=\mathbf1_{[0,\infty)}$. Write $F_\varphi$ for the distribution
function of $X$. Since $\chi_0=\mathbf1_{[0,2/3]}$ and
$\chi_1=\mathbf1_{[1/3,1]}$, for $\ell\ge1$ and with $h:=3^{-\ell}$,
\begin{equation}\label{eq:kernwindow}
 (W_3^\ell\chi_0)(x)=3^\ell\!\!\int_{x-\frac23h}^{x}\!\!\varphi_\ell,
 \qquad
 (W_3^\ell\chi_1)(x)=3^\ell\!\!\int_{x-h}^{x-\frac13h}\!\!\varphi_\ell .
\end{equation}

\begin{lemma}[The truncated densities]\label{lem:trunc}
Let $\ell\ge1$ and $c_\ell:=3^{-(\ell-1)}$.
\begin{enumerate}
\item[(a)] $\varphi_\ell(v)=\tfrac32F_{\ell-1}(3v)$ for $v\le\tfrac13$;
  in particular $\varphi_\ell$ is nondecreasing on
  $(-\infty,\tfrac13]$.
\item[(b)] $F_{\ell-1}(3x)>F_{\ell-1}(3x-c_\ell)$ for every
  $0<x<\tfrac13$.
\item[(c)] $\varphi\le\tfrac32$, and $F_\varphi(s)<\tfrac32
  \min(s,\tfrac23)$ and $F_\varphi(s)>\tfrac32\max(0,s-\tfrac13)$ for
  every $s\in(0,1)$.
\end{enumerate}
\end{lemma}

\begin{proof}
(a) Decomposing $X_\ell=\tfrac23U_1+\tfrac13X'_{\ell-1}$ with
$X'_{\ell-1}$ an independent copy and conditioning on it,
$\varphi_\ell(v)=\tfrac32\bigl(F_{\ell-1}(3v)-F_{\ell-1}(3v-2)\bigr)$.
For $v\le\tfrac13$ one has $3v-2\le-1<0$, so the second term vanishes.

(b) For $\ell=1$, $F_0=\mathbf1_{[0,\infty)}$ and $c_1=1$, so the two
sides are $1$ and $0$. For $\ell\ge2$ it suffices that $X_{\ell-1}$ have
positive density on the interior of $[0,1-c_\ell]$ and that
$(3x-c_\ell,3x]$ meet that interior in a set of positive length. The
first holds by induction on $m=\ell-1$: $\varphi_1=\tfrac32
\mathbf1_{[0,2/3]}$ is positive on $(0,1-3^{-1})$, and if $\varphi_m>0$
on $(0,1-3^{-m})$ then, by the display in~(a) without the restriction on
$v$, $\varphi_{m+1}(v)>0$ exactly when $(3v-2,3v]$ meets
$(0,1-3^{-m})$, which it does for every $0<v<1-3^{-(m+1)}$. The second
holds because $0<3x<1$: writing $\lambda:=\max(0,3x-c_\ell)$ and
$\Lambda:=\min(3x,1-c_\ell)$, if $3x\le1-c_\ell$ then
$\Lambda=3x>\lambda$, and otherwise $\Lambda=1-c_\ell>3x-c_\ell$ and
$\Lambda>0$.

(c) Two preliminaries. First, replacing each $U_j$ by $1-U_j$ sends $X$
to $1-X$, so $X$ is symmetric about $\tfrac12$ and
$F_\varphi(s)=1-F_\varphi(1-s)$. Second, $0<F_\varphi(s)<1$ for every
$0<s<1$.
Indeed, given such an $s$, choose $N$ with
$\sum_{j>N}2\cdot3^{-j}=3^{-N}<s/2$. On the event that $U_j<s/2$ for
every $j\le N$, of probability $(s/2)^N>0$, the head satisfies
$\sum_{j\le N}2U_j3^{-j}<s\sum_{j\ge1}3^{-j}=s/2$, so $X<s$ and
$F_\varphi(s)\ge\Pr(X<s)>0$; and $\Pr(X>s)=\Pr(X<1-s)>0$ by the
symmetry, so $F_\varphi(s)<1$.

Now the decomposition $X=\tfrac23U_1+\tfrac13X'$ gives
$\varphi(v)=\tfrac32\bigl(F_\varphi(3v)-F_\varphi(3v-2)\bigr)\le\tfrac32$,
and for $0<v<\tfrac13$ it gives $\varphi(v)=\tfrac32F_\varphi(3v)
<\tfrac32$, since $3v<1$. Hence $F_\varphi(s)=\int_0^s\varphi<\tfrac32s$
for every $s>0$, and $F_\varphi(s)<1$ for $s<1$, which together are the
first inequality in (c). For the second, the symmetry gives
$F_\varphi(s)=1-F_\varphi(1-s)>1-\tfrac32(1-s)=\tfrac32\bigl(s-\tfrac13\bigr)$
for $s\in(0,1)$, and $F_\varphi(s)>0$ was just shown.
\end{proof}

\begin{lemma}[Layer-cake comparison]\label{lem:layercake}
Let $\psi:[0,1]\to[0,\infty)$ be nonincreasing and let $\nu_1,\nu_2$ be
finite Borel measures on $[0,1]$, both absolutely continuous, with
distribution functions $F_i(s)=\nu_i([0,s])$. If $F_1(s)\ge F_2(s)$ for
every $s\in[0,1]$, then $\int_0^1\psi\,d\nu_1\ge\int_0^1\psi\,d\nu_2$.
If moreover $F_1(s)>F_2(s)$ for every $s\in(0,1)$ and $\psi(u)>\psi(v)$
for some $0<u<v<1$, the inequality is strict.
\end{lemma}

\begin{proof}
Write $\psi(w)=\int_0^\infty\mathbf1\{\psi(w)>\lambda\}\,d\lambda$ and
apply Tonelli, so that
$\int\psi\,d\nu_i=\int_0^\infty\nu_i\bigl(\{\psi>\lambda\}\bigr)\,d\lambda$.
Since $\psi$ is nonincreasing, $\{w\in[0,1]:\psi(w)>\lambda\}$ is, when
nonempty, an interval containing $0$ with right endpoint
$s(\lambda):=\sup\{w:\psi(w)>\lambda\}$, so it is $[0,s(\lambda)]$ or
$[0,s(\lambda))$; the two have the same $\nu_i$-measure because $\nu_i$
is absolutely continuous, and that measure is $F_i(s(\lambda))$. When
the set is empty, put $s(\lambda):=0$, and $F_i(0)=0$ for the same
reason. Hence
$\int\psi\,d\nu_i=\int_0^\infty F_i\bigl(s(\lambda)\bigr)\,d\lambda$,
and the first claim holds pointwise in $\lambda$, the hypothesis at
$s=0$ and $s=1$ covering the two degenerate values. For the second,
take $u<v$ as given. For $\lambda\in\bigl(\psi(v),\psi(u)\bigr)$, an
interval of positive length, $\psi(u)>\lambda$ forces $s(\lambda)\ge u$
and $\psi(v)\le\lambda$ forces $s(\lambda)\le v$, so $s(\lambda)$ lies
in $[u,v]\subset(0,1)$ and the two integrands differ strictly there.
\end{proof}

\subsection*{The six steps}

\emph{(7.3), $K_\ell(x,0)=K_\ell(x+3^{-(\ell+1)},1)$.} For $\ell\ge1$,
immediate from \eqref{eq:kernwindow}: replacing $x$ by $x+\tfrac13h$ in
the second integral turns its limits into those of the first. At
$\ell=0$ it is $\chi_0(x)=\chi_1(x+\tfrac13)$, which is the definition
of $\chi_1$.

\emph{(7.4), $K_\ell(x,0)>\varphi(x)>K_\ell(x,1)$ for $0<x<\tfrac13$.}
Take $\ell\ge1$; at $\ell=0$ the three quantities are $\tfrac32$,
$\varphi(x)$ and $0$. Put $\psi(v):=\varphi_\ell(x-hv)$ for
$v\in[0,1]$, so that \eqref{eq:kernwindow} reads
$K_\ell(x,0)=\tfrac32\int_0^{2/3}\psi$ and
$K_\ell(x,1)=\tfrac32\int_{1/3}^{1}\psi$, while conditioning
$X=X_\ell+hX'$ on an independent copy $X'\in[0,1]$ gives
$\varphi(x)=\int_0^1\psi(v)\varphi(v)\,dv$. All three are averages of
$\varphi_\ell$ over the window $[x-h,x]$, which lies in
$(-\infty,\tfrac13]$ because $x<\tfrac13$, so $\psi$ is nonincreasing by
Lemma~\ref{lem:trunc}(a). It also decreases strictly at an interior
pair. For $\ell=1$, $\varphi_1=\tfrac32\mathbf1_{[0,2/3]}$ and $h=\tfrac13$
give $\psi=\tfrac32$ on $[0,3x]$ and $\psi=0$ beyond, and $3x\in(0,1)$,
so any $0<u<3x<v<1$ will do. For $\ell\ge2$, $F_{\ell-1}$ is continuous,
hence so are $\varphi_\ell$ and $\psi$, and
Lemma~\ref{lem:trunc}(a) and~(b) give $\psi(0)>\psi(1)$; since
$\psi(u)\to\psi(0)$ as $u\downarrow0$ and $\psi(v)\to\psi(1)$ as
$v\uparrow1$, taking $u$ and $v$ close enough to the endpoints keeps
$\psi(u)>\psi(v)$. Now apply
Lemma~\ref{lem:layercake} twice; the two measures below have equal
total mass $1$, so the hypothesis also holds at $s=1$. With $\nu_1$ the measure of density
$\tfrac32\mathbf1_{[0,2/3]}$ and $\nu_2$ that of density $\varphi$, the
hypothesis $F_1>F_2$ on $(0,1)$ is the first inequality of
Lemma~\ref{lem:trunc}(c), and the conclusion is
$K_\ell(x,0)>\varphi(x)$. With $\nu_1$ of density $\varphi$ and $\nu_2$
of density $\tfrac32\mathbf1_{[1/3,1]}$, the hypothesis is the second
inequality there, and the conclusion is $\varphi(x)>K_\ell(x,1)$.

\emph{The unlabelled Taylor step
$\varphi(x-h)\ge\varphi(x+h)-2h\varphi'(x+h)$} needs $\varphi$ convex
near $0$, which holds because differentiating
$\varphi'(x)=\tfrac92\varphi(3x)$ gives
$\varphi''(x)=\tfrac{27}{2}\varphi'(3x)\ge0$ while $3x$ stays in the
region where $\varphi$ is nondecreasing.

\emph{(7.12), the asymptotic solution property
$\lambda\varphi_0'(t)/\bigl(a\varphi_0(at)\bigr)\to1$.} Both Wirsching
and Berg and Kr\"uppel state this without argument, and it carries
Proposition~\ref{prop:wirschingreq} and both routes to
Theorem~\ref{thm:star4}. It follows from \eqref{eq:phi0}. Write
$L:=-\log t\to\infty$, so that \eqref{eq:phi0} becomes
$\log\varphi_0(t)=C-\gamma L+\delta_{\mathrm{BK}}\log L-\beta(L+\log L)^2$
with $C=\varepsilon\log(2\beta)-\tfrac12\log(2\pi)$, and
\[
 \log\varphi_0'(t)=\log\varphi_0(t)+L+\log\Bigl(\gamma
 -\frac{\delta_{\mathrm{BK}}}{L}+2\beta(L+\log L)\bigl(1+\tfrac1L\bigr)\Bigr).
\]
Subtract the first display at $3t$, where $L$ becomes $L-c$. The linear
term contributes $-\gamma c$; the logarithmic one
$\delta_{\mathrm{BK}}\log\frac{L}{L-c}=O(1/L)$; the bracket of the
second display contributes $\log(2\beta)+\log L+O(\log L/L)$. For the
quadratic one put $A=L+\log L$ and $B=L-c+\log(L-c)$: then
$A-B=c+c/L+O(L^{-2})$ and $A+B=2L+2\log L-c+O(1/L)$, so
$A^2-B^2=2cL+2c\log L-c^2+2c+O(\log L/L)$ and, since $2\beta c=1$ and
$\beta c^2=c/2$, that term contributes $-L-\log L+\tfrac c2-1$. The $L$
and $\log L$ pieces cancel against $L$ and $\log L$ above, leaving
\begin{equation}\label{eq:w712}
 \log\frac{\varphi_0'(t)}{\varphi_0(3t)}
 \longrightarrow -\log c-\gamma c+\tfrac c2-1
 \;=\;\log\tfrac92 ,
\end{equation}
the last equality by substituting $\gamma=-2\beta-\delta_{\mathrm{BK}}
-\tfrac12$ and $\beta=1/(2c)$. At $\lambda=2/3$, $a=3$ this is exactly
(7.12). The
limit depends on $t$ only through $L$, so it is uniform over any family
of $t$ with $-\log t\to\infty$ uniformly, which is what
$\widetilde A_{\delta_5}$ supplies.

\emph{(7.13)}, introduced with ``a somewhat lengthy calculation shows'',
follows from \eqref{eq:phi0} by differentiation. With $v=\log t$,
$\varphi_0'(t)/\varphi_0(t)=B(v)/t$ where
$B(v)=\gamma+\delta_{\mathrm{BK}}/v-2\beta\bigl(v-\log(-v)\bigr)
\bigl(1-\tfrac1v\bigr)=-2\beta v\,(1+O(\log|v|/|v|))$. For
$(x_\ell)\in\widetilde A_\delta$, $x_\ell=\ell\,3^{-\ell}
\bigl(1+O(\delta/\sqrt\ell)\bigr)$ and $v=-\ell c+O(\log\ell)$, so
$\tfrac2{3^{\ell+1}}\varphi_0'(x_\ell)/\varphi_0(x_\ell)
=\tfrac23\cdot\tfrac{2\beta\ell c}{\ell}\bigl(1+o(1)\bigr)\to\tfrac23$
using $2\beta c=1$. Every error term depends on $\delta$ alone, so the
limit is uniform on the class, which is what (7.13) asserts.

\emph{(7.14), the class lifting.} Given $(x_\ell)\in\widetilde A_\delta$
and $\delta_5>\delta$, the sequences $\hat x_\ell:=x_\ell^+$ and
$\hat y_{\ell-1}:=3x_\ell^+$ lie in $\widetilde A_{\delta_5}$ for all but
finitely many $\ell$. Indeed $\lfloor3^\ell x_\ell^+\rfloor
=\lfloor3^\ell x_\ell+\tfrac13\rfloor\in\{k_\ell,k_\ell+1\}$ where
$k_\ell=\lfloor3^\ell x_\ell\rfloor$, and
$\lfloor3^{\ell-1}(3x_\ell^+)\rfloor$ is that same integer, so the window
conditions at index $\ell$ and at index $\ell-1$ both follow from
$|k_\ell-\ell|\le\delta\sqrt\ell$ once
$\delta\sqrt\ell+2\le\delta_5\sqrt{\ell-1}$, which holds for all large
$\ell$ because $\delta_5>\delta$.

\emph{One attribution.} Wirsching gives the penultimate equality of the
calculation in the proof of Theorem~\ref{thm:star4} to his (7.14) and
(7.7). His (7.7) is an identity for
$\varphi$, and the line it is applied to carries $\varphi_0$; what
licenses that step is (7.12), which contributes a factor $1+o(1)$ rather
than an equality. None of the six points needs a hypothesis the source
does not already have, and none changes $\mu=1/3$.

\bibliographystyle{plain}

\end{document}